%Version: A.Degtyarev 05-09-97  11:16:53am
%Format: plain
%
%
%	On the Pontrjagin-Viro form
%
%       by Alexander Degtyarev
%
%
%
%
%
\input amstex
\documentstyle{amsppt}
\input epsf
%\input degt.def
%\input label.def
%\nologo
%\def\today{}
%\preprint{\eightrm\today\hfil }

%%%%%%%%%% File label.def
%Version: A.Degtyarev 05-06-97  08:24:40pm
%Format: amsppt
\def\stydate{May 7, 1997}
\expandafter\ifx\csname amsppt.sty\endcsname\relax\input amsppt.sty \fi

\chardef\tempatcatcode\the\catcode`\@
\catcode`\@=11

\W@{This is LABEL.DEF by A.Degtyarev <\stydate>}
\ifx\labelloaded@\undefined\else
	\catcode`\@\tempatcatcode\let\tempatcatcode\undefined
  \message{[already loaded]}\endinput\fi
\let\labelloaded@\tempatcatcode\let\tempatcatcode\undefined
\def\labelmesg@ {LABEL.DEF: }
{\edef\temp{\the\everyjob\W@{\labelmesg@<\stydate>}}
\global\everyjob\expandafter{\temp}}

\def\stylefile@#1{\expandafter\ifx\csname#1\endcsname\relax
	\else\message{[already loaded]}\endinput\fi
	\expandafter\edef\csname#1\endcsname{\catcode`\noexpand\@\tempcat
		\toks@{}\toks@@{}\expandafter\let\csname#1\endcsname\noexpand\empty}%
	\let\tempcat\undefined\uppercase{\def\styname{#1}}%
	{\edef\temp{\the\everyjob\noexpand\W@{\styname: <\stydate>}}
		\global\everyjob\expandafter{\temp}}}
\def\loadstyle#1{\edef\next{#1}%
	\DN@##1.##2\@nil{\if\notempty{##2}\else\def\next{##1.sty}\fi}%
	\expandafter\next@\next.\@nil
	\expandafter\ifx\csname\next\endcsname\relax\input\next\fi}

\def\make@letter{\edef\t@mpcat{\catcode`\@\the\catcode`\@}\catcode`\@=11 }
{\let\head\relax\let\specialhead\relax\let\subhead\relax
\let\subsubhead\relax\let\proclaim\relax
\gdef\let@relax{\let\head\relax\let\specialhead\relax\let\subhead\relax
	\let\subsubhead\relax\let\proclaim\relax}}

\def\arabicnum#1{\number#1}

\def\Romannum#1{\expandafter\uppercase\expandafter{\romannumeral #1}}
\def\alphnum#1{\ifcase#1\or a\or b\or c\or d\else\@ialph{#1}\fi}
\def\@ialph#1{\ifcase#1\or \or \or \or \or e\or f\or g\or h\or i\or j\or
	k\or l\or m\or n\or o\or p\or q\or r\or s\or t\or u\or v\or w\or x\or y\or
	z\else\fi}
\def\Alphnum#1{\ifcase#1\or A\or B\or C\or D\else\@Ialph{#1}\fi}
\def\@Ialph#1{\ifcase#1\or \or \or \or \or E\or F\or G\or H\or I\or J\or
	K\or L\or M\or N\or O\or P\or Q\or R\or S\or T\or U\or V\or W\or X\or Y\or
	Z\else\fi}

\def\@car#1#2\@nil{#1}
\def\@cdr#1#2\@nil{#2}
\newskip\@savsk
% We add a tiny skip in order to make \@esphack work correct
\let\@ignorespaces\ignorespaces
\def\@ignorespacesp{\ifhmode
%$$  \ifdim\lastskip>\z@\else\nobreak\hskip1sp minus1sp\fi\fi\@ignorespaces}
  \ifdim\lastskip>\z@\else\penalty\@M\hskip-1sp
		\penalty\@M\hskip1sp \fi\fi\@ignorespaces}
\def\ignorespaces{\protect\@ignorespacesp}
% However, we have to redefine some control sequences
%\def~{\unskip\nobreak\ \@ignorespaces}
\def\@bsphack{\relax\ifmmode\else\@savsk\lastskip
  \ifhmode\edef\@sf{\spacefactor\the\spacefactor}\fi\fi}
\def\@esphack{\relax
  \ifx\penalty@\penalty\else\penalty\@M\fi   % if this is after \nobreak
%$$  \ifmmode\else\ifhmode\@sf{}\ifdim\@savsk>\z@\@ignorespaces\fi\fi\fi}
  \ifmmode\else\ifhmode\@sf{}\ifdim\@savsk>\z@\@ignorespacesp\fi\fi\fi}
\newread\@inputcheck
\def\@input#1{\openin\@inputcheck #1 \ifeof\@inputcheck \W@
  {No file `#1'.}\else\closein\@inputcheck \relax\input #1 \fi}

\def\eat@bs{\expandafter\eat@\string}
\def\eat@ii#1#2{}
\def\eat@iii#1#2#3{}
\def\eat@iv#1#2#3#4{}
\def\@xname#1{\expandafter\noexpand\csname\eat@bs#1\endcsname}
\def\@name#1{\csname\eat@bs#1\endcsname}
\def\@Name#1{\@ifundefined{#1}{}{\@name{#1}}}
\def\@defname#1{\expandafter\def\csname\eat@bs#1\endcsname}
\def\@gdefname{\global\@defname}
\def\@edefname#1{\expandafter\edef\csname\eat@bs#1\endcsname}
\def\@xdefname{\global\@edefname}
\long\def\@ifundefined#1#2#3{\expandafter\ifx\csname\eat@bs#1\endcsname\relax
  #2\else#3\fi}
\def\@@addto#1#2{{\toks@\expandafter{#1#2}\xdef#1{\the\toks@}}}
\def\@@addparm#1#2{{\toks@\expandafter{#1{##1}#2}%
	\edef#1{\gdef\noexpand#1####1{\the\toks@}}#1}}
\def\ST@P{@step}
\def\ST@LE{style}
\def\N@M{no}
%	#1 is the name of the counter to be defined,
%	#2 is the counter this one depends upon,
\outer\def\newcounter{\checkbrack@{\expandafter\newcounter@\@txtopt@{{}}}}
{
\gdef\newcounter@#1#2#3{{
	\toks@@\expandafter{\csname\eat@bs#2\N@M\endcsname}%
	\DN@{\alloc@0\count\countdef\insc@unt}%
	\ifx\@txtopt@\identity@\expandafter\next@\the\toks@@
		\else\if\notempty{#1}
			\global\expandafter\let\csname\eat@bs#2\N@M\endcsname#1\fi\fi
  \@ifundefined{\the\eat@bs#3}{\toks@{}}{%
		\toks@\expandafter{\csname the\eat@bs#3\endcsname.}}%
	\@xdefname{\the\eat@bs#2}{\the\toks@\noexpand\arabicnum\the\toks@@}%
  \@xdefname{#2\ST@P}{}%
  \@ifundefined{#3\ST@P}{}
  	{\edef\next@{\noexpand\@@addto\@xname{#3\ST@P}{%
			 \global\@xname{#2\N@M}\z@\@xname{#2\ST@P}}}\next@}%
	\expandafter\@@addto\expandafter\closeaux@\expandafter
		{\expandafter\\\the\toks@@}}}}
\outer\def\copycounter#1#2{%
	\@xdefname{#1\N@M}{\@xname{#2\N@M}}%
	\@xdefname{#1\ST@P}{\@xname{#2\ST@P}}%
	\@xdefname{\the\eat@bs#1}{\@xname{\the\eat@bs#2}}}
\outer\def\everystep{\checkstar@\everystep@}
\def\everystep@#1{\ifx\@numopt@\identity@\let\next@\@@addto
  \else\let\next@\gdef\fi\expandafter\next@\csname\eat@bs#1\ST@P\endcsname}
%	#1 is the counter whose style is to be changed,
\def\counterstyle#1{\@gdefname{\the\eat@bs#1}}
\def\advancecounter#1#2{\@name{#1\ST@P}\global\advance\@name{#1\N@M}#2}
\def\setcounter#1#2{\@name{#1\ST@P}\global\@name{#1\N@M}#2}
\def\counter#1{\refstepcounter#1\printcounter#1}
\def\printcounter#1{\@name{\the\eat@bs#1}}
\def\refcounter#1{\xdef\@lastmark{\printcounter#1}}
\def\stepcounter#1{\advancecounter#1\@ne}
\def\refstepcounter#1{\advancecounter#1\@ne\refcounter#1}
\def\savecounter#1{\@edefname{#1@sav}{%
	\global\@name{#1\N@M}\the\@name{#1\N@M}}}
\def\restorecounter#1{\@ifundefined{#1@sav}{}{\@name{#1@sav}}}

\def\warning#1#2{\W@{Warning: #1 on input line #2}}
\def\warning@#1{\warning{#1}{\the\inputlineno}}
\def\wrn@Protect#1#2{\warning@{\string\Protect\string#1\space ignored}}
\def\wrn@label#1#2{\warning{label `#1' multiply defined}{#2}}
\def\wrn@ref#1#2{\warning@{label `#1' undefined}}
\def\wrn@cite#1#2{\warning@{citation `#1' undefined}}
\def\wrn@command#1#2{\warning@{Preamble command \string#1\space ignored}}
\def\wrn@option#1#2{\warning@{Option \string#1\string#2\space is not supported}}
%% Disable wornings: works with
%% \Protect, \label, \ref, \cite, \command, \option,
%% (with FONT.DEF) \font
\def\nowarning#1{\@ifundefined{\wrn@\eat@bs#1}{\wrn@option\nowarning#1}
	{\expandafter\let\csname wrn@\eat@bs#1\endcsname\eat@ii}}

\def\bftext{\ifmmode\fam\bffam\else\bf\fi}
\let\@lastmark\empty
\let\@lastlabel\empty
\def\lastmark{\@lastmark}
\let\lastlabel\empty
\def\newlabel#1#2#3{{\@ifundefined{\r@-#1}{}{\wrn@label{#1}{#3}}%
  {\let\protect\noexpand\@xdefname{\r@-#1}{#2}}}}
\def\w@@@xref#1#2{%
	\@ifundefined{\r@-#1}{{\bftext??}#2{#1}{}}%
  {\expandafter\expandafter\expandafter\@car\csname r@-#1\endcsname\@nil}}%\null}}
\def\@@@xref#1{\w@@@xref{#1}\wrn@ref}
\def\@xref#1{\rom{\@@@xref{#1}}}
\let\xref\@xref
\def\pageref#1{%
	\@ifundefined{\r@-#1}{{\bftext??}\wrn@ref{#1}{}}%
  {\expandafter\expandafter\expandafter\@cdr\csname r@-#1\endcsname\@nil}}%\null}}
\def\thepage{\ifnum\pageno<\z@\romannumeral-\pageno\else\number\pageno\fi}
\def\label@#1#2{\@bsphack\edef\@lastlabel{#1}{\let\thepage\relax
  \def\protect{\noexpand\noexpand\noexpand}%
  \edef\@tempa{\write\@auxout{\string
    \newlabel{#2}{{\@lastmark}{\thepage}}{\the\inputlineno}}}%
  \expandafter}\@tempa\@esphack}
\def\label#1{\label@{#1}{#1}}
\def\@gobble{\relaxnext@
 	\DN@{\ifx[\next\DN@[####1]{}\else
 		\ifx"\next\DN@"####1"{}\else\DN@{}\fi\fi\next@}%
 	\FN@\next@}
\def\@gobblefn#1{\ifx#1[\expandafter\@gobblebr\else
  \ifx#1"\expandafter\expandafter\expandafter\@gobblequ\fi\fi}
\def\@gobblebr#1]#2{}
\def\@gobblequ#1"#2{}
{\catcode`\~\active\lccode`\~=`\@
\lowercase{\global\let\save@at=~ \gdef\protect@at{\def~{\protect\save@at}}}}
\def\Protect@@#1{\def#1{\protect#1}}
\def\disable@special{\let\@writeaux\eat@ii\let\label\eat@
	\def\footnotemark{\protect\@gobble}%
  \let\footnotetext\@gobblefn\let\footnote\@gobblefn
  \Protect@@\\\let\ifvmode\iffalse
	\let\refcounter\eat@\let\advancecounter\eat@ii\let\setcounter\eat@ii}
\def\@writeaux#1#2{\@bsphack{\disable@special\protect@at
	\def\chapter{\protect\chapter@toc}\let\thepage\relax
	\def\protect{\noexpand\noexpand\noexpand}%
	\Protect@@~\Protect@@\@@@xref\Protect@@\pageref\Protect@@\nofrills
  \edef\@tempa{\@ifundefined#1{}{\write#1{#2}}}\expandafter}\@tempa\@esphack}
\def\writeauxline#1#2#3{\@writeaux\@auxout
  {\string\@auxline{#1}{#2}{#3}{\thepage}}}
%\newtoks\thetoc@@
{\let\newwrite\relax
\gdef\@openin#1{\make@letter\@input{\jobname.#1}\t@mpcat}
\gdef\@openout#1{\global\expandafter\newwrite\csname tf@-#1\endcsname
   \immediate\openout\csname tf@-#1\endcsname \jobname.#1\relax}}
%\long\def\@writefile#1#2{\@ifundefined{\tf@-#1}{}{{\def~{\string~}%
%  \def\protect##1{\string##1 }\let\nofrills\relax
%	\immediate\write\csname tf@-#1\endcsname{#2}}}}
\def\auxlinedef#1{\@defname{\do@-#1}}
\def\@auxline#1{\@ifundefined{\do@-#1}{\expandafter\eat@iii}%
	{\expandafter\expandafter\csname do@-#1\endcsname}}
\def\begin@write#1{\bgroup\def\do##1{\catcode`##1=12 }\dospecials\do\~\do\@
	\catcode`\{=\@ne\catcode`\}=\tw@\immediate\write\csname#1\endcsname}
\def\end@writetoc#1#2#3{{\string\tocline{#1}{#2\string\page{#3}}}\egroup}
\def\do@tocline#1{%
%%	The file version
	\@ifundefined{\tf@-#1}{\expandafter\eat@iii}
		{\begin@write{tf@-#1}\expandafter\end@writetoc}
%%	The \toks version
%		\@ifundefined{\the#1@@}{}
%			{\global\addto{\the#1@}{\tocline{#2}{#3\page{#4}}}}%
}
\auxlinedef{toc}{\do@tocline{toc}}

\let\protect\empty
\def\Protect#1{\@ifundefined{#1@P@}{\PROTECT#1}{\wrn@Protect#1{}}}
\def\PROTECT#1{%
	\expandafter\let\csname\eat@bs#1@P@\endcsname#1%
	\edef#1{\noexpand\protect\@xname{#1@P@}}}
\def\pdef#1{\edef#1{\noexpand\protect\@xname{#1@P@}}\@defname{#1@P@}}

\Protect\operatorname
\Protect\operatornamewithlimits
\Protect\qopname@
\Protect\qopnamewl@
\Protect\text
\Protect\topsmash
\Protect\botsmash
\Protect\smash
\Protect\widetilde
\Protect\widehat
\Protect\thetag
\Protect\therosteritem
% Fonts:
\Protect\Cal
\Protect\Bbb
\Protect\bold
\Protect\slanted
\Protect\roman
\Protect\italic
\Protect\boldkey
\Protect\boldsymbol
\Protect\frak
\Protect\goth
\Protect\dots
% Symbols
\Protect\cong
\Protect\lbrace \let\{\lbrace
\Protect\rbrace \let\}\rbrace
\let\root@P@@\root \def\root@P@#1{\root@P@@#1\of}
\def\root#1\of{\protect\root@P@{#1}}

\let\@frills@\identity@
\let\@txtopt@\identyty@
\let\@numopt@\identity@
\def\frills{\ignorespaces\@txtopt@}
\def\numberline{\@numopt@}
\newif\if@write\@writetrue
\def\checkstar@#1{\DN@{\@writetrue
  \ifx\next*\DN@####1{\let\@numopt@\eat@\checkstar@@{#1}}%
	  \else\DN@{\let\@numopt@\identity@#1}\fi\next@}\FN@\next@}
\def\checkstar@@#1{\DN@{%
  \ifx\next*\DN@####1{\@writefalse#1}%
	  \else\DN@{\@writetrue#1}\fi\next@}\FN@\next@}
\def\checkfrills@#1{\DN@{%
  \ifx\next\nofrills\DN@####1{#1}\def\@frills@####1{####1\nofrills}%
	  \else\DN@{#1}\let\@frills@\identity@\fi\next@}\FN@\next@}
\def\checkbrack@#1{\DN@{%
	\ifx\next[\DN@[####1]{\def\@txtopt@########1{####1}%
%% This is nice for the `1.1. Theorem.' style, but not for `Theorem 1.1.'
%		\if\notempty{####1}\else\def\@frills@########1{########1\nofrills}\fi
		#1}%
	\else\DN@{\let\@txtopt@\identity@#1}\fi\next@}\FN@\next@}
\def\check@therstyle#1#2{{\DN@{#1}\ifx\@txtopt@\identity@\else
		\DNii@##1\@therstyle{}\def\@therstyle{\DN@{#2}\nextii@}%
		\def\nextiii@##1##2\@therstyle{\expandafter\nextii@##1##2\@therstyle}%
    \expandafter\nextiii@\@txtopt@\@therstyle.\@therstyle\fi
	\expandafter}\next@}

\newif\if@theorem
\let\@therstyle\eat@
\def\@headtext@#1#2{{\disable@special\let\protect\noexpand
	\def\chapter{\protect\chapter@rh}\Protect@@\nofrills
	\edef\@tempa{\noexpand\@frills@\noexpand#1{#2}}\expandafter}\@tempa}
\let\AmSrighthead@\rightheadtext
\def\rightheadtext{\checkfrills@{\@headtext@\AmSrighthead@}}
\let\AmSlefthead@\leftheadtext
\def\leftheadtext{\checkfrills@{\@headtext@\AmSlefthead@}}
% #1 refers to the style,
% #2 refers to the style in the toc,
% #3 refers to the counter,
% #4 represents the AmSTeX's counterpart (we use \end... because
%   of \outer-ness), and
% #5 is the text to be typeset
\def\@head@@#1#2#3#4#5{\@Name{\pre\eat@bs#1}\if@theorem\else
	\@frills@{\csname\expandafter\eat@iv\string#4\endcsname}\relax
    	\ifx\protect\empty\@name{#1font@}\fi\fi
  	\@name{#1\ST@LE}{\counter#3}{\ignorespaces#5\unskip}%
  \if@write\writeauxline{toc}{\eat@bs#1}{#2{\counter#3}{#5}}\fi
	\if@theorem\else\expandafter#4\fi
	\ifx#4\endhead\ifx\@txtopt@\identity@\else
		\headmark{\@name{#1\ST@LE}{\counter#3}{\frills{}}}\fi\fi
	\@Name{\post\eat@bs#1}\ignorespaces}
\ifx\undefined\endhead\Invalid@\endhead\fi
\def\@head@#1{\checkstar@{\checkfrills@{\checkbrack@{\@head@@#1}}}}
% #1 is the name,
% #2 is the counter, and
% #3 is the title text
\def\@thm@@#1#2#3{\@Name{\pre\eat@bs#1}%
	\@frills@{\csname\expandafter\eat@iv\string#3\endcsname}
  	{\@theoremtrue\check@therstyle{\@name{#1\ST@LE}}\frills
			{\counter#2}\@theoremfalse}%
	\expandafter\envir@stack\csname end\eat@bs#1\endcsname
	\@name{#1font@}\@Name{\post\eat@bs#1}\ignorespaces}
\def\@thm@#1{\checkstar@{\checkfrills@{\checkbrack@{\@thm@@#1}}}}
% #1 is the name,
% #2 refers is the counter,
% #3 is the name of the corresponding table (lof, lot, etc.)
% #4 is either \topcaption or \botcaption
%	#5 is the caption text
\def\@capt@@#1#2#3#4#5\endcaption{\bgroup
	\edef\@tempb{\global\footmarkcount@\the\footmarkcount@
    \global\@name{#2\N@M}\the\@name{#2\N@M}}%
	\def\shortcaption##1{\def\sh@rtt@xt####1{##1}}\let\sh@rtt@xt\identity@
	\DN@##1##2##3{\false@\fi\iftrue}%
	\ifx\@frills@\identity@\else\let\notempty\next@\fi
  #4{\@tempb\@name{#1\ST@LE}{\counter#2}}\@name{#1font@}#5\endcaption
  \if@write\writeauxline{#3}{\eat@bs#1}{{} \@name{#1\ST@LE}{\counter#2}%
    \if\notempty{#5}.\enspace\fi\sh@rtt@xt{#5}}\fi\egroup}
\def\@capt@#1{\checkstar@{\checkfrills@{\checkbrack@{\@capt@@#1}}}}
\let\captiontextfont@\empty

\ifx\undefined\subsubheadfont@\def\subsubheadfont@{\it}\fi
\ifx\undefined\proclaimfont\def\proclaimfont{\sl}\fi
\ifx\undefined\proclaimfont@\let\proclaimfont@\proclaimfont\fi
\def\proclaimfont{\proclaimfont@}
\ifx\undefined\definitionfont@\def\AmSdeffont@{\rm}
	\else\let\AmSdeffont@\definitionfont@\fi
\ifx\undefined\remarkfont@\def\remarkfont@{\rm}\fi

\def\newfont@def#1#2{\@ifundefined{#1font@}
	{\@xdefname{#1font@}{\@xname{.\expandafter\eat@iv\string#2font@}}}{}}
% #1 is the name (and the style),
% #2 is the style in the toc,
% #3 is the counter, and
% #4 is the AmSTeX's counterpart to be used (in the \end... form):
\def\newhead@#1#2#3#4{{%
	\gdef#1{\@therstyle\@therstyle\@head@{#1#2#3#4}}\newfont@def#1#4%
	\@ifundefined{#1\ST@LE}{\@gdefname{#1\ST@LE}{\headstyle}}{}%
	\@ifundefined{#2}{\gdef#2{\headtocstyle}}{}%
  \@@addto\moretocdefs@{\\#1#1#4}}}
\outer\def\newhead#1{\checkbrack@{\expandafter\newhead@\expandafter
	#1\@txtopt@\headtocstyle}}
% #1 is the default title (like Theorem, Lemma, etc.),
% #2 is the name to be defined,
% #3 refers to the counter (which should be defined separately),
% #4 is the AmSTeX's counterpart: \endproclaim, \endremark, \endAmSdef
\outer\def\newtheorem#1#2#3#4{{%
	\gdef#2{\@thm@{#2#3#4}}\newfont@def#2#4%
	\@xdefname{\end\eat@bs#2}{\noexpand\revert@envir
		\@xname{\end\eat@bs#2}\noexpand#4}%
  \@ifundefined{#2\ST@LE}{\@gdefname{#2\ST@LE}{\proclaimstyle{#1}}}{}}}%
% #1 is the default title (like Figure, Table, etc.),
% #2 is the name to be defined,
% #3 refers to the counter (which should be defined separately),
% #4 is the name of the corresponding table (toc, lof, lot, etc.)
% #5 is either \topcaption or \botcaption
\outer\def\newcaption#1#2#3#4#5{{\let#2\relax
  \edef\@tempa{\gdef#2####1\csname end\eat@bs#2\endcsname}%
	\@tempa{\@capt@{#2#3{#4}#5}##1\endcaption}\newfont@def#2\endcaptiontext%
  \@ifundefined{#2\ST@LE}{\@gdefname{#2\ST@LE}{\captionstyle{#1}}}{}%
  \@@addto\moretocdefs@{\\#2#2\endcaption}\newtoc{#4}}}
{
% #1 is the name of the table (toc, lof, lot, etc.)
\outer\gdef\newtoc#1{{%
	\expandafter\ifx\csname do@-#1\endcsname\relax
%%	The \toks version
%    \global\expandafter\newtoks\csname the#4@@\endcsname
    \global\auxlinedef{#1}{\do@tocline{#1}}{}%
    \@@addto\tocsections@{\make@toc{#1}{}}\fi}}}

\toks@\expandafter{\itembox@}
\toks@@{{\let\therosteritem\identity@\let\rm\empty
  \edef\next@{\edef\noexpand\@lastmark{\therosteritem@}}\expandafter}\next@}
\edef\itembox@{\the\toks@@\the\toks@}
\def\firstitem@false{\let\iffirstitem@\iffalse
	\global\let\lastlabel\@lastlabel}

\def\rosteritemref#1{\hbox{\therosteritem{\@@@xref{#1}}}}
\def\local#1{\label@\@lastlabel{\lastlabel-i#1}}
\def\loccit#1{\rosteritemref{\lastlabel-i#1}}
\def\xRef@P@{\gdef\lastlabel}
\def\xRef#1{\@xref{#1}\protect\xRef@P@{#1}}

\def\iref@P@{\gdef\lastref}
\def\itemref#1#2{\rosteritemref{#1-i#2}\protect\iref@P@{#1}}
\def\iref#1{\@xref{#1}\itemref{#1}}
\def\ditto#1{\rosteritemref{\lastref-i#1}}

\def\eqtag{\tag\counter\equation}
\def\eqref#1{\thetag{\@@@xref{#1}}}
\def\tagform@#1{\ifmmode\hbox{\rm\else\rom{\fi
	(\ignorespaces#1\unskip)\iftrue}\else}\fi}

\let\AmSfnote@\makefootnote@
\def\makefootnote@#1{{\let\footmarkform@\identity@
  \edef\next@{\edef\noexpand\@lastmark{#1}}\expandafter}\next@
  \AmSfnote@{#1}}

\def\clearpage{\ifnum\insertpenalties>0\line{}\fi\vfill\supereject}

\def\find@#1\in#2{\let\found@\false@
	\DNii@{\ifx\next\@nil\let\next\eat@\else\let\next\nextiv@\fi\next}%
	\edef\nextiii@{#1}\def\nextiv@##1,{%
    \edef\next{##1}\ifx\nextiii@\next\let\found@\true@\fi\FN@\nextii@}%
	\expandafter\nextiv@#2,\@nil}
\def\include@#1{{\ifx\@auxout\@subaux\DN@{\errmessage
		{\labelmesg@ Only one level of \string\include\space is supported}}%
	\else\edef\@tempb{#1}\@numopt@\clearpage
		\xdef\@include@{\relax\bgroup\ifx\@numopt@\identity@\clearpage
			\else\let\noexpand\@immediate\noexpand\empty\fi}%
	  \DN@##1 {\if\notempty{##1}\edef\@tempb{##1}\DN@####1\eat@ {}\fi\next@}%
  	\DNii@##1.{\edef\@tempa{##1}\DN@####1\eat@.{}\next@}%
		\expandafter\next@\@tempb\eat@{} \eat@{} %
  	\expandafter\nextii@\@tempb.\eat@.%
		\let\next@\empty
	  \if\expandafter\notempty\expandafter{\@tempa}%
		  \edef\nextii@{\write\@mainaux{%
  			\noexpand\string\noexpand\@input{\@tempa.aux}}}\nextii@
  		\@ifundefined\@includelist{\let\found@\true@}
				{\find@\@tempa\in\@includelist}%
			\if\found@\@ifundefined\@noincllist{\let\found@\false@}
				{\find@\@tempb\in\@noincllist}\else\let\found@\true@\fi
			\if\found@\@ifundefined{\close@#1}{}{\DN@{\csname close@#1\endcsname}}%
			\else\xdef\@auxname{\@tempa}\xdef\@inputname{\@tempb}%
				\@numopt@{\Wcount@@{open@\@auxname}}%
				\global\let\@auxout\@subaux
	 	  	\@numopt@\immediate\openout\@auxout=\@tempa.aux
				\@numopt@\immediate\write\@auxout{\relax}%
 	  		\DN@{\let\end\endinput\@input\@inputname \let\end\endmain@
          \@include@\closeaux@@\egroup\make@auxmain}\fi\fi\fi
  \expandafter}\next@}
\def\include{\checkstar@\include@}
\def\includeonly#1{\edef\@includelist{#1}}
\def\noinclude#1{\edef\@noincllist{#1}}

\newwrite\@mainaux
\newwrite\@subaux
\def\make@auxmain{\global\let\@auxout\@mainaux
  \xdef\@auxname{\jobname}\xdef\@inputname{\jobname}}
\begingroup
\catcode`\(\the\catcode`\{\catcode`\{=12
\catcode`\)\the\catcode`\}\catcode`\}=12
%\gdef\closeaux@@((%
%	\def\\##1(\@Waux(\global##1=\the##1))%\edef\@tempb(\closeaux@)%
%	\edef\@tempa(\@Waux(%
%		\gdef\expandafter\string\csname cl@\@auxname\endcsname{)%
%		\\\pageno\\\footmarkcount@\closeaux@\@Waux(})\@immediate\closeout\@auxout)%
%  \expandafter)\@tempa)
\gdef\Wcount@@#1((%
	\def\\##1(\@Waux(\global##1=\the##1))%
	\edef\@tempa(\@Waux(%
		\string\expandafter\gdef\string\csname\space#1\string\endcsname{)%
		\\\pageno\\\footmarkcount@\closeaux@\@Waux(}))\expandafter)\@tempa)
\endgroup
\def\closeaux@@{\Wcount@@{close@\@auxname}\@immediate\closeout\@auxout}
\let\closeaux@\empty
\let\@immediate\immediate
\def\@Waux{\@immediate\write\@auxout}
\def\readaux{\checkbrack@\readaux@}
\def\readaux@{%
	\W@{>>> \labelmesg@ Run this file twice to get x-references right}%
	\global\everypar{}%
	{\def\\##1{\global\let##1\relax}%
		\def\/##1{\gdef##1{\wrn@command##1{}}}%
		\disablepreambule@cs}%
	\make@auxmain\make@letter{\setboxz@h{\@input{\@txtopt@{\jobname.aux}}}}\t@mpcat
  \immediate\openout\@auxout=\jobname.aux%
	\immediate\write\@auxout{\relax}\global\let\end\endmain@}
\everypar{\global\everypar{}\readaux}
{\toks@\expandafter{\topmatter}
\global\edef\topmatter{\noexpand\readaux\the\toks@}}
\let\@@end@@\end

\def\endmain@{\clearpage\@immediate\closeout\@auxout
	\make@letter\def\newlabel##1##2##3{}\@input{\jobname.aux}%
 	\W@{>>> \labelmesg@ Run this file twice to get x-references right}%
 	\@@end@@}
\def\disablepreambule@cs{\\\disablepreambule@cs}

\def\proof{\checkfrills@{\checkbrack@{%
	\check@therstyle{\@frills@{\demo}{\frills{Proof}}{}}
		{\frills{}\envir@stack\endremark\envir@stack\enddemo}%
  \envir@stack\endproof\ignorespaces}}}
\def\endproof{\nofrillscheck{\frills@{\qed}\revert@envir\endproof\enddemo}}

\let\AmSref\ref
\let\AmSrefstyle\refstyle
\let\plaincite\cite
\def\citei@#1,{\citeii@#1\eat@,}
\def\citeii@#1\eat@{\w@@@xref{#1}\wrn@cite}
\def\cite#1{\protect\plaincite{\citei@#1\eat@,\unskip}}
\def\refstyle#1{\AmSrefstyle{#1}\uppercase{%
	\ifx#1A\relax \def\@ref@##1{\AmSref\xdef\@lastmark{##1}\key##1}%
  	\else\ifx#1C\relax \def\@ref@{\AmSref\no\counter\refno}%
		\else\def\@ref@{\AmSref}\fi\fi}}
\refstyle A
\newcounter\refno\null
\gdef\Refs{\checkstar@{\checkbrack@{\csname AmSRefs\endcsname
  \nofrills{\frills{References}%
  \if@write\writeauxline{toc}{vartocline}{\frills{References}}\fi}%
  \def\ref{\@ref@}\ignorespaces}}}
\let\ref\xref

\def\tocsections@{\make@toc{toc}{}}
\let\moretocdefs@\empty
\def\newtocline@#1#2#3{%
  \edef#1{\def\@xname{#2line}####1{\expandafter\noexpand
      \csname\expandafter\eat@iv\string#3\endcsname####1\noexpand#3}}%
  \@edefname{\no\eat@bs#1}{\let\@xname{#2line}\noexpand\eat@}%
	\@name{\no\eat@bs#1}}
\def\maketoc#1#2{\Err@{\Invalid@@\string\maketoc}}
\def\newtocline#1#2#3{\Err@{\Invalid@@\string\newtocline}}
\def\make@toc#1#2{\penaltyandskip@{-200}\aboveheadskip
	\if\notempty{#2}
		\centerline{\headfont@\ignorespaces#2\unskip}\nobreak
  	\vskip\belowheadskip \fi
%%	The file version
	\@openin{#1}\@openout{#1}%
%%	The \toks version
%	\@ifundefined{\the#1@@}{}
%		{\the\@name{\the#1@@}\global\@name{\the#1@@}}{}%
	\vskip\z@}
\def\contents{\readaux\checkfrills@{\checkbrack@{\@contents@}}}
\def\@contents@{\toc@{\frills{Contents}}\envir@stack\endcontents%
	\def\nopagenumbers{\let\page\eat@}\let\newtocline\newtocline@
  \def\tocline##1{\csname##1line\endcsname}%
  \edef\caption##1\endcaption{\expandafter\noexpand
    \csname head\endcsname##1\noexpand\endhead}%
	\ifmonograph@\def\vartoclineline{\Chapterline}%
		\else\def\vartoclineline##1{\sectionline{{} ##1}}\fi
  \let\\\newtocline@\moretocdefs@
	\ifx\@frills@\identity@\def\\##1##2##3{##1}\moretocdefs@
		\else\let\tocsections@\relax\fi
	\def\\{\unskip\space\ignorespaces}\let\maketoc\make@toc}
\def\endcontents{\tocsections@\vskip-\lastskip\revert@envir\endcontents
	\endtoc}

% \selectf@nt is for future extensions (like RUSSIAN.TEX)
\@ifundefined\selectf@nt{\let\selectf@nt\identity@}{}
\def\textonlyfont@#1#2{%
	\@defname{#1@P@}{\RIfM@\Err@{Use \string#1\space only in text}%
		\else\edef\f@ntsh@pe{\string#1}\selectf@nt#2\fi}%
	\edef#1{\noexpand\protect\@xname{#1@P@}\empty}}
\tenpoint

% #1 is the name of the switch to be defined
%	#2 is the default font switch
\def\newshapeswitch#1#2{\gdef#1{\selectsh@pe#1#2}\PROTECT#1}
% #1 is the name of the switch
% #2 is the current shape
% #3 is the shape to be used
\def\shapeswitch#1#2#3{\@gdefname{#1\string#2}{#3}}
% These shapes are used by \rom
\shapeswitch\rm\bf\bf  \shapeswitch\rm\tt\tt  \shapeswitch\rm\smc\smc
\newshapeswitch\em\it
% These shapes are used by \em and \emph
\shapeswitch\em\it\rm  \shapeswitch\em\sl\rm
\def\selectsh@pe#1#2{\relax\@ifundefined{#1\f@ntsh@pe}{#2}
	{\csname\eat@bs#1\f@ntsh@pe\endcsname}}

\def\@itcorr@{\leavevmode\skip@\lastskip\unskip\/%
  \ifdim\skip@=\z@\else\hskip\skip@\fi}
\def\rom@P@#1{\@itcorr@{\selectsh@pe\rm\rm#1}}
\def\rom{\protect\rom@P@}
%%%%	\Rom will unconditionally switch to \rm, like \rom in AmSppt
\def\Rom@P@#1{\@itcorr@{\rm#1}}
\def\Rom{\protect\Rom@P@}
{\catcode`\-=11
\gdef\wr@index#1{\@writeaux\tf@-idx{\string\indexentry{#1}{\thepage}}}
\gdef\wr@glossary#1{\@writeaux\tf@-glo{\string\glossaryentry{#1}{\thepage}}}}
\let\index\wr@index
\let\glossary\wr@glossary
%% This version would also write to .idx
%\def\emph{\@itcorr@\bgroup\em\checkstar@\emph@}
%\def\emph@#1{\@numopt@{\wr@index{#1}}\ignorespaces#1\unskip\egroup
%  \DN@{\DN@{}\ifx\next.\else\ifx\next,\else\DN@{\/}\fi\fi\next@}\FN@\next@}
\def\emph#1{\@itcorr@\bgroup\em\ignorespaces#1\unskip\egroup
  \DN@{\DN@{}\ifx\next.\else\ifx\next,\else\DN@{\/}\fi\fi\next@}\FN@\next@}
\def\makequoteactive{\catcode`\"\active}
{\makequoteactive\gdef"{\FN@\quote@}
\gdef\quote@{\ifx"\next\DN@"##1""{\wr@glossary{##1}}\else
	\DN@##1"{\wr@index{##1}}\fi\next@}}
\def\MakeIndex{\@openout{idx}}
\def\MakeGlossary{\@openout{glo}}

%%%%%%%%%%%%%%%%%%%%%%	Just helpful things	%%%%%%%%%%%%%%%%%%%%%%%%%%%%
\def\endofpar#1{\ifmmode\ifinner\endofpar@{#1}\else\eqno{#1}\fi
	\else\leavevmode\endofpar@{#1}\fi}
\def\endofpar@#1{\unskip\penalty\z@\null\hfil\hbox{#1}\hfilneg\penalty\@M}

\newdimen\normalparindent\normalparindent\parindent
\def\firstparindent#1{\everypar\expandafter{\the\everypar
  \parindent\normalparindent\everypar{}}\parindent#1\relax}

%% Commands to disable
\@@addto\disablepreambule@cs{%
	\\\readaux
	\/\Monograph
	\/\MakeIndex
	\/\MakeGlossary
}

\catcode`\@\labelloaded@

%%%%%%%%%%%%%%%%%%%%%%%% Definitions for my papers %%%%%%%%%%%%%%%%%%%%%%%

\def\headstyle#1#2{\numberline{#1.\enspace}#2}
\def\headtocstyle#1#2{\numberline{#1.}\space #2}

\def\specialtocstyle#1#2{#2}
\newcounter\section\null
\newcounter\subsection\section
\newcounter\subsubsection\subsection
\newhead\specialsection[\specialtocstyle]\null\endspecialhead
\newhead\section\section\endhead
\newhead\subsection\subsection\endsubhead
\newhead\subsubsection\subsubsection\endsubsubhead
\def\firstappendix{\global\sectionno=0 %
  \counterstyle\section{\Alphnum\sectionno}%
	\global\let\firstappendix\empty}

\def\appendixtocstyle#1#2{\space\numberline{Appendix #1.\enspace}#2}
\newhead\appendix[\appendixtocstyle]\section\endhead

\let\endAmSdef\enddefinition
\def\proclaimstyle#1#2{\numberline{#2.\enspace}\frills{#1}}
\copycounter\thm\subsubsection
%\newcounter[\subsubsectionno]\thm\subsection
\theorem\thm\endproclaim
\proposition\thm\endproclaim
\lemma\thm\endproclaim
\corollary\thm\endproclaim
\definition\thm\endAmSdef
\example\thm\endAmSdef

\def\captionstyle#1#2{\frills{\numberline{#1 #2}}}
\newcounter\figure\null
\newcounter\table\null
\newcaption{Figure}\figure\figure{lof}\botcaption
\newcaption{Table}\table\table{lot}\topcaption

\copycounter\equation\subsubsection
%\newcounter[\subsubsectionno]\equation\subsection

%%%%%%%%%% End of label.def

\let\ge\geqslant
\let\le\leqslant
\def\C{{\Bbb C}}
\def\R{{\Bbb R}}
\def\Z{{\Bbb Z}}
\def\Q{{\Bbb Q}}

\def\Cp#1{\C{\roman p}^{#1}}
\def\Rp#1{\R{\roman p}^{#1}}
{\catcode`\@=11
\gdef\Hom{\qopname@{Hom}}
\gdef\Tors{\qopname@{Tors}}
\gdef\Im{\qopname@{Im}}			%	I use these
\gdef\Ker{\qopname@{Ker}}
\gdef\Fix{\qopname@{Fix}}
\gdef\tr{\qopname@{tr}}
\gdef\inj{\qopname@{in}}
\gdef\id{\qopname@{id}}
\gdef\pr{\qopname@{pr}}
\gdef\rel{\qopname@{rel}}
\gdef\conj{\qopname@{conj}}
\global\let\sminus\smallsetminus
\gdef\emptyset{\varnothing}
}

\let\cedilla\c

\def\+{\mathbin{\scriptstyle\sqcup}}
\def\scup{\mathbin{\scriptstyle\cup}}
\def\scap{\mathbin{\scriptstyle\cap}}
\def\stimes{\mathbin{\scriptstyle\times}}
\def\star{\mathbin*}
\let\<\langle
\let\>\rangle
\def\bv{\operatorname{bv}}
\let\vr\bv

\def\bX{\bar X}
\def\binj{\operatorname{\bar{\text{\rm\i n}}}}

\def\case#1#2{\par{\it Case #1\/{\rm:}\enspace#2}}

\def\rH{\RH{r@!@!}}
\def\rd{\Rd{r@!}}
\def\RH#1{{}^{#1}\!@!H}
\def\Rd#1{{}^{#1}\!@!d}
\def\iH{\RH\infty}
\def\rB{\RB{r@!@!}}
\def\RB#1{{}^{#1}\!@!B}
\def\iB{\RB\infty}
\def\rZ{\RZ{r@!@!}}
\def\RZ#1{{}^{#1}\!@!Z}
\def\iZ{\RZ\infty}
\def\rE{\RE{r@!@!}}
\def\RE#1{{}^{#1}\!@!E}

\def\Sq{\operatorname{Sq}}

\def\CF{\Cal F}
\def\hidepar#1{\vcenter{\hbox{$\scriptscriptstyle#1$}}}
\def\CFF#1{\CF_{\smash{\hidepar[#1\hidepar]}}}
\def\CP{\Cal P}
\def\CPP#1{\CP_{\smash{\hidepar[#1\hidepar]}}}
\def\gM{\frak M}
\def\gN{\frak N}
\def\gS{\frak S}
\def\gC{\frak C}
\def\gA{\frak A}
\def\gtl{\frak l}
\let\Gk\varkappa
\let\Gd\delta
\let\Ge\varepsilon
\let\Gl\lambda
\let\Gm\mu
\def\Br{\operatorname{Br}}
\def\Gr{\operatorname{Gr}}
\def\CW{{\sl CW}}
\def\Card{\operatorname{Card}}
\def\sint{\operatorname{int}}
\def\Spin{\operatorname{Spin}}
\def\Pin{\operatorname{Pin}}
\def\gm{\operatorname{\frak{q}}}
\def\ind{\operatorname{ind}}

\def\lk{\operatorname{lk}}

\def\er{E_{\R}\futurelet\nexts\getsup}
\def\xr{X_{\R}\futurelet\nexts\getsup}

\def\c{\futurelet\nexts\conjs}
\def\index#1{^{\botsmash{(#1)}}}
\def\conjs{\ifx\nexts1\def\nexts##1{t\index1}\else
 \ifx\nexts2\def\nexts##1{t\index2}\else
 \def\nexts{\mathop{\roman{conj}}}\fi\fi\nexts}
\def\eri{\er\index i}
\def\xri{\xr\index i}
\def\ci{t\index i}
\def\getsup{\ifx\nexts1\def\nexts##1{\index1}\else
 \ifx\nexts2\def\nexts##1{\index2}\else\let\nexts\relax\fi\fi\nexts}

\def\QR#1#2{Q_{#2}\index{#1}}
\def\qr#1#2{\frak{q}_{#2}\index{#1}}

\def\tX{X}
\def\tY{Y}
\def\tZ{Z}
\def\bY{\bar Y}
\def\bZ{\bar Z}
\def\bP{\bar P}
\def\bQ{\bar Q}

\def\ttX{\smash{\tilde X}}

\def\I#1{\roman{I}_{#1}}
\def\II{\roman{II}}
\def\Iu{\I{u}}

\def\Pic{\operatorname{Pic}}

\Protect\er

\topmatter

\title     %   Title
On the Pontrjagin-Viro form
\endtitle
%\keftheadtext{} %   Short title

\author
Alexander Degtyarev
\endauthor

\address
Steklov Mathematical Institute, %\newline\indent
St.~Petersburg, Russia and \newline\indent
Bilkent University, Ankara, Turkey
\endaddress

\email
degt\,\@\,fen.bilkent.edu.tr, degt\,\@\,pdmi.ras.ru
\endemail

%\date    %   Date
%
%\enddate

%\dedicatory    %   Dedicatory
%
%\enddedicatory



\thanks    %   Thanks
An essential part of this work was done during my visits to
\emph{Universit\'e Louis Pasteur}, Strasbourg, \emph{Universit\'a di
Trento}, and \emph{Universit\'e de Rennes \rom{I}}.
\endthanks

\abstract    %   Abstract
A new invariant, the Pontrjagin-Viro form, of algebraic surfaces is
introduced and studied. It is related to various
Rokhlin-Guillou-Marin forms and generalizes Mikhalkin's complex
separation. The form is calculated for all real Enriques surfaces for
which it is well defined.
\endabstract

\endtopmatter

\document

\section*{Introduction}

In this paper we introduce a new invariant, so called
\emph{Pontrjagin-Viro form}, of a real algebraic surface or, more
generally, a closed smooth $4n$-manifold~$X$ with
involution~$c\:X\to X$.  The invariant, which is only well defined in
certain special cases, is a quadratic function $\CP\:\CF\to\Z/4$,
where $\CF\subset H_*(\Fix c;\Z/2)$ is a subgroup of the total
homology of the fixed point set of~$c$ (or, in the case of
an algebraic surface, of the real part of the surface). We mainly
concentrate on the case $\dim X=4$; in this case $\CP$ turns out to
be closely related to the \emph{Rokhlin-Guillou-Marin forms}
(see~\ref{2.2}) of various characteristic surfaces in~$X$ and $X/c$
and, thus, is a direct generalization of the notion of \emph{complex
separation} introduced by G.~Mikhalkin~\cite{Mikhalkin}. (Mikhalkin's
complex separation is defined when $H_1(X;\Z/2)=0$.) The relation to
the Rokhlin-Guillou-Marin form gives a number of congruences which
the Pontrjagin-Viro form must satisfy (see~\ref{4.2}).

This work was mainly inspired by our study of real Enriques surfaces
(joint work with I.~Itenberg and V.~Kharlamov). Recall that an
\emph{Enriques surface} is a complex analytic surface~$E$ with
$\pi_1(E)=\Z/2$ and $2c_1(E)=0$. Such a surface is called \emph{real}
if it is supplied with an anti-holomorphic involution $\c\:E\to E$;
the fixed point set $\er=\Fix\c$ is called the \emph{real part}
of~$E$. The set of components of the real part of a real Enriques
surface naturally splits into two disjoint \emph{halves} $\er1$,
$\er2$ (see~\ref{5.1}); this splitting is a deformation invariant of
pair $(E;\c)$.

The topology of real Enriques surfaces is studied in~\cite{DK1}
and~\cite{DK2}, where they are classified up to homeomorphism of the
triad $(\er;\er1,\er2)$. Currently, we know the classification up to
deformation equivalence (which is the strongest equivalence relation
from the topological point of view); a preliminary report is found
in~\cite{DK3}; details will appear in~\cite{DIK}. For a technical
reason real Enriques surfaces are divided in~\cite{DK3} into three
types, \emph{hyperbolic}, \emph{parabolic}, and \emph{elliptic},
according to whether the minimal Euler characteristic of the
components of~$\er$ is negative, zero, or positive, respectively. It
turns out that in most cases a real Enriques surface is determined up
to deformation equivalence by such classical invariants as the
homeomorphism type of the triad \smash{$(\er;\er1,\er2)$} and whether
the fundamental classes~$[\er]$ and~$[\eri]$, $i=1,2$, vanish or are
characteristic in the homology of~$E$ or some auxiliary manifolds.
The few exceptions, mainly $M$-surfaces of parabolic and elliptic
types, differ by the Pontrjagin-Viro form.

In this paper the Pontrjagin-Viro form is calculated for all real
Enriques surfaces for which it is well-defined (see
Section~\ref{s7}). There is a necessary condition
($\chi(\er)=0\bmod8$) and certain sufficient conditions
(Lemma~\ref{sufficient}) for $\CP$ to be well defined, and, when
defined, $\CP$ must satisfy certain congruences
(Proposition~\ref{E-RGM}) which follow from the general congruences
in~\ref{4.2}. The result of the calculation can be roughly stated as
follows (see Theorems~\ref{PV-ep}, \ref{PV-hyp}, and~\ref{th-other}
for the precise statements): {\proclaimfont Consider a triad
$(\er;\er1,\er2)$ with $\chi(\er)=0\bmod8$. Any \rom(partial\rom)
quadratic form $\CP\:H_*(\er1)\oplus H_*(\er2)\to\Z/4$ satisfying the
congruences of Proposition~\ref{E-RGM} can be realized as the
Pontrjagin-Viro form of a real Enriques surface. If $(\er;\er1,\er2)$
does\/ {\bf not} satisfy the sufficient conditions of
Lemma~\ref{sufficient}, it can also be realized by a real Enriques
surface not admitting Pontrjagin-Viro form.}
Note that, in fact, the deformation type of a surface admitting
Pontrjagin-Viro form is determined by the topology of
$(\er;\er1,\er2)$ and the isomorphism type of
$\CP\:H_*(\er1)\oplus H_*(\er2)\to\Z/4$ (see~\cite{DIK}).

Originally in order to distinguish nonequivalent real Enriques
surfaces we calculated the Pontrjagin-Viro form by explicitly
constructing membranes in $E/\c$. In Section~\ref{s6} I develop a
different approach, which facilitates the calculation and, on the
other hand, covers all real Enriques surfaces which are of interest.
The approach is applicable to a specific construction (which, as is
shown in~\cite{DIK}, produces all $M$-surfaces of elliptic and
parabolic types): the surface in question is constructed starting
from a pair of real curves $P$, $Q$ on a real rational surface~$\tZ$,
and the Pontrjagin-Viro form is given in terms of the topology of
their real parts $(\tZ_\R;P_\R,Q_\R)$. The fact that $\CP$ is related
to the complex orientation of the branch curve was indicated to me by
G.~Mikhalkin.

\subsection*{Contents}
Section~\ref{s1} introduces the primary tool, so called Kalinin's
spectral sequence. Section~\ref{s2} reminds the basic notions related
to quadratic forms and Rokhlin-Guillou-Marin form of a characteristic
surface. The Pontrjagin-Viro form is introduced in Section~\ref{s3},
and its basic properties, including the congruences, are studied in
Section~\ref{s4}. In Section~\ref{s5} the general results are
transferred to real Enriques surfaces. In Section~\ref{s6} we
calculate the Pontrjagin-Viro form of a real Enriques surface
obtained by a specific construction, using so called Donaldson's
trick; these results are applied in Section~\ref{s7} to produce the
complete list of possible values of the Pontrjagin-Viro form on a
real Enriques surface.

\subsection*{Acknowledgements}
I am thankful to I.~Itenberg and V.~Kharlamov for many hours of
fruitful discussions of the subject; the results of this paper depend
essentially upon our joint work on real Enriques surfaces. I am also
thankful to G.~Mikhalkin for his helpful remarks.

\subsection*{Notation}
Unless stated otherwise, all homology and cohomology groups are with
coefficients in~$\Z/2$. We freely denote by $2\:\Z/2\to\Z/4$ the
nontrivial homomorphism, as well as the induced homomorphisms
$H_*(\,\cdot\,;\Z/2)\to H_*(\,\cdot\,;\Z/4)$ etc.

Given a vector bundle~$\xi$, we denote by~$w_i(\xi)$ and $u_i(\xi)$
the Stiefel-Whitney and Wu classes, respectively. $w=1+w_1+\dots$ and
$u=1+u_1+\dots$ are the corresponding total classes.  If $X$ is a
smooth manifold and $\tau X$ its tangent bundle, we abbreviate
$w_i(\tau X)=w_i(X)$ and $u_i(\tau X)=u_i(X)$. For a smooth
submanifold $V\subset X$ we denote by $\nu V=\nu_XV$ its normal
bundle in~$X$.

Let $X$ be a closed manifold of dimension~$n$. Then $[X]\in H_n(X)$
is its fundamental class and $\<X\>\in H_0(X)$ is the $0$-class
defined by the union of points, one in each component of~$X$.
(Certainly, the latter definition applies to any polyhedron.) The
Poincar\'e duality $\scap[X]\:H^i(X)\to H_{n-i}(X)$ is denoted by
$D_X=D$.

If $X$ is a complex manifold, $c_i(X)\in H_{2i}(X;\Z)$ stand for its
Chern classes and $K_X$, for both the canonical class in $\Pic(X)$
and its image $D_Xc_1(K_X)$ in $H_2(X)$ (so that we can write
$[D]=K_X$ for a divisor~$D$).

\section{Kalinin's spectral sequence}\label{s1}

\subsection{Basic concepts}\label{1.1}
Let $X$ be a good topological space (say, a finite dimensional
\CW-complex) and $c\:X\to X$ an involution. Unless stated otherwise,
we assume~$X$ connected. Denote by~$F$ the fixed point set $\Fix c$
and by~$\bX$, the orbit space $X/c$. Let $\pr\:X\to \bX$ be the
projection and $\inj\:F\hookrightarrow X$ and
$\binj\:F\hookrightarrow\bX$ the inclusions.

Recall that the \emph{Borel-Serre spectral sequences} $\rE_{p,q}$ and
$\rE^{p,q}$ are the Serre spectral sequences of the fibration
$S^\infty\times_cX\to\Rp\infty$, where $S^\infty\times_cX$
is the \emph{Borel construction}
$S^\infty\times X/(\bold s,x)\sim(-\bold s,cx)$. As shown
in~\cite{Kalinin}, multiplication by the generator
$h\in H^1(\Rp\infty)$ establishes isomorphisms
$\rE_{p,q+1}\to\rE_{p,q}$ and $\rE^{p,q}\to\rE^{p,q+1}$ for $p\gg0$
and thus produces stabilized spectral sequences $(\rH_*,\rd_*)$ and
$(\rH^*,\rd^*)$, which we call \emph{Kalinin's spectral
sequences} of~$(X,c)$, so that
\roster
\item\local1
$\RH1_*=H_*(X)$ and $\RH1^*=H^*(X)$,
\item\local2
$\Rd1_*=(1+c_*)$ and $\Rd1^*=(1+c^*)$,
\item\local3
$\rH_*\Rightarrow H_*(F)$ and $\rH^*\Rightarrow H^*(F)$.
\endroster
An alternative, geometrical, description of Kalinin's spectral
sequences and related objects is found in~\cite{DK2}.

The convergence in~\loccit3 means that there is an increasing
filtration $\{\CF^p\}$ on $H_*(F)$, a decreasing
filtration~$\{\CF_p\}$ on~$H^*(F)$, and homomorphisms
$\bv_p\:\CF^p\to\iH_p$ and $\bv^p\:\iH^p\to H^*(F)/\CF_{p-1}$ which
establish isomorphisms of the graded groups. (Note that in general
the filtrations do not respect the grading on $H_*(F)$ and $H^*(F)$.)
We will call~$\bv_p$ and~$\bv^p$ the \emph{Viro homomorphisms}; often
they will be considered as additive relations (partial homomorphisms)
$H_*(F)\dasharrow H_p(X)$ and $H^p(X)\dasharrow H^*(F)$.

\subsection{Multiplicative structures}\label{1.2}
The homology and cohomology versions of Ka\-li\-nin's spectral
sequences are dual to each other, i.e., $\rH^p=\Hom(\rH_p,\Z/2)$ and
$\rd^p=\Hom(\rd_p,\id_{\Z/2})$. The cup- and cap-products convert
$\rH^*$ and $\rH_*$ to a graded $\Z/2$-algebra and a graded
$\rH^*$-module, respectively, so that all the differentials
except~$\Rd1$ are differentiations. Furthermore, if $X$ is a closed
$n$-manifold and $F\ne\varnothing$, the Poincar\'e duality~$D_X$
induces isomorphisms $D\:\rH^p\to\rH_{n-p}$, and in the usual way one
can define intersection pairings
$\star\:\rH_p\otimes\rH_q\to\rH_{p+q-n}$. The induced (via~$\bv_*$)
pairing on the graded group $\Gr^*_{\CF}H_*(F)$ is called
\emph{Kalinin's intersection pairing}. The ordinary intersection
pairing on $H_*(F)$ will be denoted by~$\circ$.

\theorem[Theorem \rm(see~\cite{DK2})]\label{circ}
Let $X$ be a smooth closed $n$-manifold and~$c\:X\to X$ a smooth
involution. Then for $a\in\CF^p$ and $b\in\CF^q$ one has
$w(\nu F)\scap(a\circ b)\in \CF^{p+q-n}$
and
$$
\bv_p a\circ\bv_q b=\bv_{p+q-n}[w(\nu F)\scap(a\circ b)].
$$
\endtheorem

\subsection{Relation to the Smith exact sequence}\label{1.3}
Recall that the \emph{Smith exact sequences} of~$(X,c)$ are the exact
sequences
$$
\gather
@>>>H_{p+1}(\bX,F)@>\Delta>>H_p(\bX,F)\oplus H_p(F)
@>\tr_*>>H_p(X)@>\pr_*>>H_p(\bX,F)@>>>\rlap{\,,}\\
@>>>H^p(\bX,F)@>\pr^*>>H^p(X)@>\tr^*>>H^p(\bX,F)\oplus H^p(F)
@>\Delta>>H^{p+1}(\bX,F)@>>>\rlap{\,.}
\endgather
$$
The connecting homomorphisms~$\Delta$ are given by
$x\mapsto\omega\scap x\oplus\partial x$ (in homology) and
$x\oplus f\mapsto\omega\scup x+\delta f$ (in cohomology), where
$\omega\in H^1(\bX\sminus F)$ is the characteristic class of the
double covering $X\sminus F\to\bX\sminus F$. In~\cite{DK2} it is shown
that Kalinin's spectral sequences can be derived from the Smith exact
sequences. In this paper we only need the corresponding description
of the differentials and Viro homomorphisms:

\theorem[Theorem \rm(see~\cite{DK2})]\label{Smith}
The differentials~$\rd_p$ and~$\rd^p$, considered as additive
relations $H_p(X)\dasharrow H_{p+r-1}(X)$ and
$H^p(X)\dasharrow H^{p-r+1}(X)$, are given by
$$
\rd_*=\tr_*\circ\iota\circ(\Delta^{-1}\circ\iota)^{r-1}\circ\pr_*,\qquad
\rd^*=\pr^*\circ(\pi\circ\Delta^{-1})^{r-1}\circ\pi\circ\tr^*,
$$
where $\iota\:H_p(\bX,F)\to H_p(\bX,F)\oplus H_p(F)$ and
$\pi\:H^p(\bX,F)\oplus H^p(F)\to H^p(\bX,F)$ are, respectively, the
inclusion and the projection.

A \rom(nonhomogeneous\rom) class $x=\sum_{i\le p} x_i$,
$x_i\in H_i(F)$, belongs to~$\CF^p$ if and only if there are elements
$y_i\in H_i(\bX,F)$ such that $\Delta(y_{i+1})=y_i\oplus x_i$ for
$i<p$. In this case $\vr_px=\tr_*(y_p\oplus x_p)$
\rom(modulo the indeterminacy subgroup\rom).

A class $x\in H^p(X)$ survives to~$\iH^p$ if and only if
$\tr^*x=y^p\oplus x^p$ extends to a sequence
$y^i\oplus x^i\in H^i(\bX,F)\oplus H^i(F)$, $i\le p$, such that
$y^{i+1}=\Delta(y^i\oplus x^i)$ for $i<p$. In this case
$\bv^px=\sum_{i\le p}x^i\bmod\CF_{p-1}$.
\endtheorem

\subsection{Groups~$\rB$ and~$\rZ$}\label{1.4}
Let $\rB_p\subset\rZ_p\subset H_p(X)$ be the pull-backs of
$\Im\Rd{r-1}_p$ and $\Ker\Rd{r-1}_p$, respectively, so that
$\rH_p=\rZ_p/\rB_p$. Denote $\iB_p=\bigcup_r\rB_p$ and
$\iZ_p=\bigcap_r\rZ_p$. Then $\iH_p=\iZ_p/\iB_p$. There are
obvious cohomology analogues $\rB^p\subset\rZ^p\subset H^p(X)$, and
$\rH^p=\rZ^p/\rB^p$ for $1\le p\le\infty$.

\proposition\label{B,Z}
One has
$\iZ_p=\Ker[\pr_*\:H_p(X)\to H_p(\bX,F)]$ and
$\iB^p=\Im[\pr^*\:H^p(\bX,F)\to H^p(X)]$.
\endproposition

\proof
The statement follows from Theorem~\ref{Smith}. Since all the spaces
involved are finite dimensional, both the Smith exact sequences
terminate.  Hence, an element $x\in H_p(X)$ is annihilated by
all~$\rd_p$, $p>0$, if and only if $(\iota\circ\pr_*)(x)=0$.
Similarly, for any element $x\in H^p(\bX,F)$ the multiple image
$\Delta^{r-1}(x)$ belongs to $\Im(\pi\circ\tr^*)$ for $r\gg0$.
\endproof

\corollary\label{B,Z'}
If $X$ is a closed smooth manifold, $c$ is a smooth involution, and
$F\ne\emptyset$, then
$\iB_*=\Im[\tr_*\:H_*(\bX\sminus F)\to H_*(X)]$ and
$\iZ^*=\Ker[\tr^*\:H^*(X)\to H^*(\bX\sminus F)]$.
\endcorollary

\proof
This follows from~\ref{B,Z} and Poincar\'e duality.
\endproof

\subsection{Miscellaneous statements}\label{1.5}
In this section we state several simple results needed in the sequel.

\proposition[Proposition \rm(see~\cite{DK2})]\label{Bockstein}
Denote by~$\Sq_1$ the homology Bockstein homomorphism. Then
for any class $x=\sum_{i\le p}x_i\in\CF^p$, $x_i\in H_i(F)$, one has
$$
\tsize
\Sq_1\bv_px=
\bv_{p-1}\big(\Sq_1x+\sum_iix_i\big)=
\bv_{p-1}\big(\Sq_1x+\sum_i(i+1)x_i\big).
$$
\rom(In particular, the classes in parentheses belong
to~$\CF^{p-1}$.\rom)
\endproposition

\proposition\label{char.class}
Let~$X$ be an oriented closed smooth $n$-manifold, $H_1(X)=0$, and
$B\subset X$ a $c$-invariant oriented closed smooth submanifold of
pure codimension~$2$ such that $[B]=0$ in $H_{n-2}(X)$. Assume that
$c$ reverses the co-orientation of~$B$.  Let, further, $p:Y\to X$ be
the \rom(unique\rom) double covering of~$X$ branched over~$B$ and
$\omega_p$ its characteristic class. Then for $x\in H_1(F\sminus B)$
one has $\<\omega_p,x\>=\vr_2x\circ\frac12[B]$ \rom(where
$\frac12[B]$ is obtained by dividing by~$2$ the integral class
$[B]\in H_{n-2}(X;\Z)$\,\rom).
\endproposition

\proof
Realize~$x$ by an oriented simple loop~$\gtl$. After multiplying it by
an odd integer we may assume that $\gtl$ bounds an oriented
membrane~$\gM$ in~$X$, which may be chosen transversal to~$B$. Then
$\<\omega_p,x\>=\Card(\gM\cap B)\bmod2$. On the other hand,
$[\gM\cup c(\gM)]=\bv_2x$ and the statement follows from
$\Card(\gM\cap B)=\frac12[\gM\cup c(\gM)]\circ[B]$. (Note that $c$
reverses the orientation of~$\gM$.)
\endproof

\proposition\label{sphere}
Assume that $X$ is a closed $4$-manifold, $H_1(X;\Z)=0$, and $F$ is a
surface. Then $\bX$ is a $\Z$-homology $4$-sphere if and only if
$c$ is an $M$-involution \rom(i.e., $\dim H_*(F)=\dim H_*(X)$\,\rom)
and $F$ is connected. If this is the case, $c_*$ acts as
multiplication by~$(-1)$ on $H_2(X;\Z)$.
\endproposition

\proof
Assume that $F\ne\varnothing$ (as otherwise $H_1(\bX)=\Z/2$). Then
$H_1(\bX;\Z)=0$ and the first statement follows from comparing the
Euler characteristics using the Riemann-Hurwitz formula. For the
second statement observe that $H_*(\bX;\Q)$ is the
$c_*$-skew-invariant part of $H_*(X;\Q)$; this determines the action
of~$c_*$ on $H_*(X;\Z)\subset H_*(X;\Q)$.
\endproof

\section{Rokhlin-Guillou-Marin congruence}\label{s2}

\subsection{Quadratic forms and Brown invariant}\label{2.1}
The results of this section should be found in most textbooks in
arithmetics; see also~\cite{vdBlij}, \cite{Brown}, \cite{GM},
or~\cite{KV}.

Let~$V$ be a $\Z/2$-vector space and $\circ\:V\otimes V\to\Z/2$ a
symmetric bilinear form. A function $q\:V\to\Z/4$ is called a
\emph{quadratic extension} of~$\circ$ if
$q(x+y)=q(x)+q(y)+2(x\circ y)$ for all $x,y\in V$. The pair $(V,q)$ is
called a \emph{quadratic space}. (Obviously, $\circ$ is recovered
from~$q$.) A quadratic space is called \emph{nonsingular} if the
bilinear form is nonsingular, i.e., $V^\perp=0$; it is called
\emph{informative} if $q|_{V^\perp}=0$. The following is
straightforward:

\proposition\label{q-space}
Let $V$ be a $\Z/2$-vector space and $\circ\:V\otimes V\to\Z/2$ a
symmetric bilinear form. Then
\roster
\item\local1
$q(x)\equiv x^2\bmod2$ for any $x\in V$ and any quadratic
extension~$q$ of~$\circ$\rom;
\item\local2
quadratic extensions of~$\circ$ form an affine space over
$V\spcheck=\Hom(V,\Z/2)$ via $(q+l)(x)=q(x)+2l(x)$ for
$l\in V\spcheck$ and $x\in V$\rom;
\item\local3
if $\circ$ is nonsingular, its quadratic extensions form
an affine space over~$V$ via $(q+v)(x)=q(x)+2(v\circ x)$ for
$v,x\in V$.
\endroster
\endproposition

The \emph{Brown invariant} $\Br(V,q)$ (or just $\Br q$) of a
nonsingular quadratic space is the $(\!\!{}\bmod8)$-residue defined
by
$$
\tsize
\exp\bigl(\tfrac14i\pi\Br q\bigr)=
2^{\frac12\dim V}\sum_{x\in V}\exp\bigl(\tfrac12i\pi q(x)\bigr).
$$
This notion extends to informative spaces: since $q$ vanishes
on~$V^\perp$, it descends to a quadratic form $q'\:V/V^\perp\to\Z/4$,
and one lets $\Br q=\Br q'$.

\proposition\label{Br}
For any informative quadratic space $(V,q)$ one has\rom:
\roster
\item\local1
$\Br q\equiv\dim(V/V^\perp)\bmod2$\rom;
\item\local2
$\Br q\equiv q(u)\bmod4$ for any characteristic element $u\in V$\rom;
\item\local3
$\Br(q+v)=\Br q-2q(v)$ for any $v\in V$
\rom(see~\iref{q-space}3\rom)\rom;
\item\local4
$\Br q=0$ if and only if $(V,q)$ is\/ \emph{null cobordant}, i.e.,
there is a subspace $H\subset V$ such that $H^\perp=H$ and $q|_H=0$.
\endroster
The Brown invariant is additive\rom: for any pair~$(V_i,q_i)$,
$i=1,2$, of quadratic spaces one has
$\Br(V_1\oplus V_2,q_1\oplus q_2)=\Br(V_1,q_1)+\Br(V_2,q_2)$.
\endproposition

\proposition\label{BrL}
Let $L$ be a unimodular integral lattice \rom(i.e., a free abelian
group with a nonsingular symmetric bilinear form
$L\otimes L\to\Z$\rom). Let $V=L\otimes\Z/2$ and define a quadratic
form $q\:V\to\Z/4$ via $q(x)=\bar x^2\bmod4$ for $x\in V$ and
$\bar x\in L$ such that $\bar x\equiv x\bmod2L$. Then
$\Br(V,q)\equiv\sigma(L)\bmod8$.
\endproposition

A subspace~$W$ of an informative quadratic space~$(V,q)$ is called
\emph{informative} if $W^\perp\subset W$ and $q|_{W^\perp}=0$.
(Clearly, an informative subspace is an informative space; hence, its
Brown invariant is well defined.)

\proposition\label{inf-space}
If $W$ is an informative subspace of an informative quadratic
space~$(V,q)$, then $\Br(W,q|_W)=\Br(V,q)$.
\endproposition

\remark{Remark}
The notion of informative subspace still makes sense if the quadratic
form~$q$ is defined only on~$W$. Proposition~\ref{inf-space} can then
be interpreted as follows: the Brown invariant of any extension
of~$q$ to a quadratic form on~$V$ equals~$\Br q$.
\endremark

\subsection{Rokhlin-Guillou-Marin congruence
\rm(see~\cite{GM})}\label{2.2}
Let $Y$ be an oriented closed smooth $4$-manifold and $U$ a
\emph{characteristic surface} in~$Y$, i.e., a smooth closed
$2$-submanifold (not necessarily orientable) with $[U]=u_2(Y)$ in
$H_2(Y)$. Denote by $i\:U\hookrightarrow Y$ the inclusion and let
$K=\Ker[i_*\:H_1(U)\to H_1(Y)]$.  Then there is a natural function
$\gm\:K\to\Z/4$, which is a quadratic extension of the intersection
index form on $H_1(U)$. We call it the \emph{Rokhlin-Guillou-Marin
form} of~$(Y,U)$. It can be defined as follows: pick a class $x\in K$
and realize it by a union~$\gtl$ of disjoint simple closed smooth
loops in~$U$.  It spans an immersed surface~$\gM$ in~$Y$, which can be
chosen normal to~$U$ along $\gtl=\partial\gM$ and transversal to~$U$
at its inner points.  (Such a surface is called a \emph{membrane}.)
Consider a normal line field~$\xi$ on~$\gtl$ tangent to~$U$ and define
the \emph{index} $\ind\gM\in\frac12\Z$ as one half of the obstruction
to extending~$\xi$ to a normal line field on~$\gM$.  (Since
$\tau\gM\oplus\nu\gM$ is an oriented vector bundle, the obstruction
is well defined as an integer. If $\gtl$ is two-sided in~$U$, the
index is usually defined using vector fields instead of line fields;
this explains the factor~$\frac12$.) Then
$\gm(x)=2\ind\gM+2\Card(\sint\gM\cap F)\bmod4$.

\theorem[Theorem \rm(see~\cite{GM})]\label{RGMcong}
Let $Y$, $U$, and $(K,\gm)$ be as above. Then $(K,\gm)$ is an
informative subspace of $H_1(U)$ and
$2\Br\gm=\sigma(Y)-U\circ U\bmod16$, where $U\circ U$ stands for the
normal Euler number of~$U$ in~$Y$.
\endtheorem

\remark{Remark}
There is an alternative construction of the
Rokhlin-Guillou-Marin form. Since $U$ is characteristic,
$Y\sminus U$ admits a $\Spin$-structure which does not extend through
any component of~$U$. Its restriction to the boundary of a tubular
neighborhood of~$U$ induces in a natural way a $\Pin^-$-structure
on~$U$ (cf.~\cite{Finashin}), which defines a quadratic form~$q$
on~$H_1(U)$. It is not difficult to see that $q$ is well defined up
to adding elements of $\Im[i^*\:H^1(Y)\to H^1(U)]$
(see~\iref{q-space}2) and, hence, its restriction to~$K$ does not
depend on the choice of a $\Spin$-structure; it coincides with~$\gm$.
\endremark

\section{Pontrjagin-Viro form}\label{s3}

\subsection{Definition of the Pontrjagin-Viro form}\label{3.1}
The \emph{Pontrjagin square} is the cohomology operation
$P^{2n}\:H^{2n}(X)\to H^{4n}(X;\Z/4)$ uniquely defined by the
following properties (see, e.g.,~\cite{Pontrjagin}):
\roster
\item\local1
$P^{2n}(x+y)=P^{2n}(x)+P^{2n}(y)+2(x\scup y)$ for any
$x,y\in H^{2n}(X)$;
\item\local2
$P^{2n}(x)\equiv x^2\bmod2$ for any $x\in H^{2n}(X)$;
\item\local3
$P^{2n}(\bar x\bmod2)=\bar x^2$ for any $\bar x\in H^{2n}(X;\Z/4)$.
\endroster
Constructively $P^{2n}$ can be defined via
$P^{2n}(x)=(\bar x\scup_0\bar x+\bar x\scup_1\delta\bar x)\bmod4$,
where $\bar x\in C^{2n}(X;\Z)$ is an integral cochain
representing~$x$ and $\scup_i$ are the cup-$i$-products used in one
of the definition of Steenrod squares.

From now on we assume that $X$ is a connected oriented closed smooth
manifold of dimension~$4n$ and $c$ is a smooth involution.
Denote by $P_{2n}\:H_n(X)\to\Z/4$ the composition
$$
H_{2n}(X)@>D_X^{-1}>>H^{2n}(X)@>P^{2n}>>H^{4n}(X;\Z/4)@>\scap[X]>>\Z/4.
$$

\proposition[Proposition-Definition]
If $P_{2n}(\iB_{2n})=0$, then $P_{2n}$ descends to a well-defined
quadratic function $\iH_{2n}\to\Z/4$. The composition of this function
and the Viro homomorphism $\bv_{2n}\:\CF^{2n}\to\iH_{2n}$ is denoted
by~$\CP$ and is called the \emph{Pontrjagin-Viro form}. It is a
quadratic extension of Kalinin's intersection form
$\star\:\CF^{2n}\otimes\CF^{2n}\to\Z/2$, i.e.,
$\CP(x+y)=\CP(x)+\CP(y)+2(x\star y)$ for any $x,y\in\CF^{2n}$.
\endproposition

\proof
The statement follows immediately from the fact that
$\iB_*\circ\iZ_*=0$ (where $\circ$ stands for the intersection
pairing) and property~\iref{3.1}1 above.
\endproof

\proposition
$\CP$ is well defined if and only if $u_{2n}(\bX\sminus F)=0$.
\endproposition

\proof
The statement is a consequence of Corollary~\ref{B,Z'} and the
obvious relation
$P_{2n}(\tr_*x)=2(x\circ x)=2\<u_{2n}(\bX\sminus F),x\>$ for
$x\in H_{2n}(\bX\sminus F)$. The latter follows from
properties~\iref{3.1}1 and~\ditto2 of Pontrjagin squares and the fact
that $\tr_*$, when restricted to the manifold $\bX\sminus F$,
coincides with the inverse Hopf homomorphism~$\pr^!$.
\endproof

\corollary\label{char.surface}
Assume that $F$ has pure dimension $\dim X-2$, so that $\bX$ is a
manifold. Then $\CP$ is well defined if and only if
$D_{\bX}u_{2n}(\bX)$ belongs to the image of
$\binj_*\:H_{2n}(F)\to H_{2n}(\bX)$.
\endcorollary

\proposition\label{congF}
If $\CP$ is well defined, $\Br\CP\equiv\sigma(X)\bmod8$.
\endproposition

\proof
By definition, $\CP$ is well defined if and only if $\iZ_{2n}$ is
an informative sub\-space of $(H_{2n}(X),P_{2n})$. (Due to the
Poincar\'e duality $\iZ_{2n}^\perp=\iB_{2n}$, see~\cite{DK2}.) Hence,
$\Br\CP=\Br P_{2n}$. On the other hand, $H_{2n}(X;\Z)\otimes\Z/2$ is
also an informative subspace, and the congruence follows from
Proposition~\ref{BrL} applied to $H_{2n}(X;\Z)/\Tors$.
\endproof

\subsection{Some sufficient conditions}\label{3.2}

\proposition\label{theta}
Let $c$ be orientation preserving. Then
$P_{2n}$ descends to~$\RH2_{2n}$ if and only if the $2n$-dimensional
component of $\inj_!(u(\tau F)u^{-1}(\nu_X F))$ equals $u_{2n}(X)$.
\endproposition

\proof
As known, the $2n$-dimensional component of
$\inj_!(u(\tau F)u^{-1}(\nu_X F))$ coincides with the characteristic
class~$\theta_{2n}$ of the twisted intersection form
$(x,y)\mapsto x\circ c_*y$ (see, e.g.,~\cite{CM}).  On the other
hand, for $x\in H_{2n}(X)$ one has
$$
P_{2n}(\Rd1_{2n}x)=P_{2n}(x+c_*x)=2P_{2n}(x)+2(x\circ c_*x)=
2\<\theta_{2n}+u_{2n},x\>
$$
(since $P_{2n}(x)\equiv x^2\bmod2$), and the statement follows.
\endproof

\corollary\label{M-case}
A necessary condition for~$\CP$ to be well defined is that
$\theta_{2n}$, the $2n$-dimensional component of
$\inj_!(u(\tau F)u^{-1}(\nu_X F))$, must coincide with $u_{2n}(X)$. The
following are sufficient conditions\rom:
\roster
\item\local1
$c$ is an $M$-involution \rom(i.e., $\rd_*=0$ for $r\ge1$\rom)\rom;
\item\local2
$\theta_{2n}=u_{2n}(X)$ and $c$ is $\Z/2$-Galois maximal
\rom(i.e., $\rd_*=0$ for $r\ge2$\rom)\rom;
\item\local3
$\theta_{2n}=u_{2n}(X)$ and $H_i(X)=0$ for $0<i<2n$.
\endroster
\endcorollary

\remark{Remark}
If $\dim F\le2n$, the condition of~\ref{theta} (and, hence, the
necessary condition of~\ref{M-case}) reduces to $[F]_{2n}=Du_{2n}(X)$
in $H_{2n}(X)$, where $[F]_{2n}$ is the fundamental class of the
union of $2n$-dimensional components of~$F$.  In this case~\ref{theta}
can be proved using the following observation:
\endremark

\proposition[Proposition \rm(V.~Arnol$'$d)]\label{Arnold}
If $\dim X=2k$ is even and $\dim F\le k$, then the fundamental
class~$[F]_k$ of the union of $k$-dimensional components of~$F$
realizes in~$H_k(X)$ the characteristic class of the twisted
intersection form $(x,y)\mapsto x\circ c_*y$.
\endproposition

\remark{Remark}
If $\dim X=4$ and $F$ is a surface, \ref{theta} follows also from the
projection formula $D_Xu_2(X)=\tr_*D_{\bX}u_2(\bX)+[F]$: since
$D_{\bX}u_2(\bX)$ comes from~$F$, its pull-back in~$X$ is zero.
\endremark

\remark{Remark}
Note that conditions~\loccit1, \loccit2 in~\ref{M-case} do not
require actual calculation of the differentials. Indeed, from
Kalinin's spectral sequence it follows that $c$ is an $M$-involution
if and only if $\dim H_*(F)=\dim H_*(X)$ (and in this case
$\theta_{2n}=u_{2n}(X)$, as $c_*=\id$ and the twisted intersection
form coincides with the ordinary one). Furthermore, $c$ is
$\Z/2$-Galois maximal if and only if $\dim H_*(F)=\dim\RH1_*$, the
latter group being equal to $\Ker(1+c_*)/\Im(1+c_*)$.
\endremark

\subsection{Membranes}\label{3.3}
The first statement, which is a direct consequence of the
definitions, calculates the Pontrjagin square in a $4$-manifold.

\lemma\label{P2}
Let $X$ be an oriented closed smooth $4$-manifold and $\gM\to X$ an
immersed closed surface. Then
$P_2[\gM]=\gM\circ\gM+2\chi(\gM)\bmod4$, where
$\gM\circ\gM=e(\mu\gM)+2i\bmod4$ is the normal Euler number of~$\gM$
plus twice the number of its self-intersection points $\!{}\bmod4$.
\rom(If $\gM$ is oriented, $\gM\circ\gM=[\gM]^2_\Z\bmod4$, where
$[\gM]_\Z\in H_2(X;\Z)$ is the integral class realized by~$\gM$.\rom)
\endlemma

Next two statements provide for geometrical means of
calculating~$\bv_2$ and, hence, the Pontrjagin-Viro form in a
$4$-manifold.

\lemma\label{bv2}
Let $\gM$ be a closed surface with involution~$c$ so that $\Fix c$
consists of several two-sided circles $l_1,\dots,l_p$, several
one-sided circles $n_1,\dots,n_q$, and several simple isolated points
$P_1,\dots,P_r$. Then $[\gM]=\bv_2\Gk$, where
$$\textstyle
\Gk=\sum[l_i]+\sum[n_j]+\sum[P_k]+\sum\<n_i\>.
$$
\endlemma

\proof
Without loss of generality we may assume that $\gM$ is connected.
Then, due to~\ref{Arnold}, $\vr_1\bigl(\sum[l_i]+\sum[n_j]\bigr)$
equals $w_1(\gM)$ in~$\iH_1$; hence,
$q\equiv\dim H_*(\gM)\bmod2$. Due to the Smith congruence
$\chi(\Fix c)\equiv\chi(\gM)\bmod2$, also
$r\equiv\dim H_*(\gM)\bmod2$.  Thus, $\vr_0\Gk=0$ and $\vr_1\Gk$ is
well defined.  Now one can easily check that $\vr_1\Gk$ annihilates
$\iH_1$ (which is generated by the images under~$\bv_1$ of~$[l_i]$,
$[n_j]$, and elements of the form $\<Q_1-Q_2\>$, $Q_1,Q_2\in\Fix c$).
Since Kalinin's intersection form is nondegenerate, $\vr_1\Gk=0$.
Hence, $\vr_2\Gk$ is well defined, and it must coincide with the only
nontrivial element $[\gM]\in\iH_2$.
\endproof

\corollary\label{Pbv2}
Let $\gM$ and~$\Gk$ be as in~\ref{bv2}.
If $\gM$ is equivariantly immersed to a topological space~$X$
with involution, then $[\gM]$ realizes~$\vr_2\Gk$. If, further,
$\dim X=4$ and $\CP$ is well defined, then
$\CP(\Gk)=\gM\circ\gM+2\chi(\gM)\bmod4$.
\endcorollary

\section{Congruences}\label{s4}

\subsection{Characteristic surfaces in~$\bX$}\label{4.1}
Let us assume that $X$ is an oriented closed smooth $4$-manifold,
$c\:X\to X$ is a smooth orientation preserving involution, and
$F=\Fix c\ne\varnothing$ has pure dimension~$2$. Under these
assumptions $\bX$ is also an oriented closed manifold.

We keep the notation introduced in Section~\ref{s1}. In addition,
denote by $\CFF{i}^p$ and $\hat\CFF{i}^p$, respectively, the
intersection $\CF^p\cap H_i(F)$ and the projection of $\CF^p$ to
$H_i(F)$. Recall that the connecting homomorphism~$\Delta$ of the
homology Smith exact sequence is given by
$y\mapsto\omega\scap y\oplus\partial y$, where
$\omega\in H^1(\bX\sminus F)$ is the characteristic class of the
covering $X\sminus F\to\bX\sminus F$. Since the covering $X\to\bX$ is
branched along~$F$, one has $\partial D_{\bX}\omega=[F]$.

\lemma\label{im0-2}
$\vr_2\CFF0^2=\tr_*H_2(\bX)\bmod\iB_2$. Furthermore, for
$y\in H_2(\bX)$ one has $\tr_*y=\bv_2(\binj^!y)\bmod\iB_2$\rom;
in particular, $\CP(\binj^!y)=2y^2\bmod4$ provided that $\CP$ is well
defined.
\endlemma

\proof
As follows from~\ref{Smith}, the image~$\vr_2\CFF0^2$ consists of all
elements of the form $\tr_*y_2$, where $y_2\in H_2(\bX,F)$ extends to
a sequence $y_i\in H_i(\bX,F)$, $i=0,1,2$, such that
$\Delta(y_2)=y_1$ and $\Delta(y_1)=y_0\oplus x_0$ for some
$x_0\in H_0(F)$. Thus, $y_1$ may be an arbitrary element, and the
only restriction to~$y_2$ is $\partial y_2=0$, i.e., $y_2=\rel y$ for
some $y\in H_2(\bX)$. For the `furthermore' part observe that
$\tr_*y=\tr_*y_2=\bv_2x_0$ and
$$
x_0=\partial(\omega\scap\rel y)=
\partial(D\omega\scap D^{-1}y)=
[F]\scap\binj^*D^{-1}y=D_F\binj^*D^{-1}y=\binj^!y,
$$
where $D=D_{\bX}$.
\endproof

\lemma\label{im0-1}
$\vr_1\CFF0^1=\tr_*H_1(\bX,F)\bmod\iB_1$.
\endlemma

\lemma\label{im1-2}
$\hat\CFF1^2=\Ker[\binj_*\:H_1(F)\to H_1(\bX)]$.
\endlemma

\proof[Proof of~\ref{im0-1} and~\ref{im1-2}]
As above, the statements follow directly from~\ref{Smith}.
\endproof

\corollary
$\Ker[\binj_*\:H_2(F)\to H_2(\bX)]$ is the annihilator of~$\CFF0^2$
with respect to the intersection index pairing
$H_2(F)\otimes H_0(F)\to\Z/2$ \rom(or, equivalently, Kalinin's
intersection pairing $H_2(F)\otimes\CFF0^2\to\Z/2$\rom).
\endcorollary

\proof
Let $K=\Ker[\binj_*\:H_2(F)\to H_2(\bX)]$. Since the restrictions of
the ordinary intersection index pairing and Kalinin's intersection
pairing to $H_2(F)\otimes\CFF0^2$ coincide (see~\ref{circ}), it
suffices to verify that $u\in K$ if and only if $\bv_2u$ annihilates
$\bv_2\CFF0^2=\tr_*H_2(\bX)$ in~$\iH_2$. For $y\in H_2(\bX)$ one has
$\inj_*u\circ\tr_*y=\binj_*u\circ y$; this product vanishes for all
$y\in H_2(\bX)$ if and only if $\binj_*u=0$.
\endproof

\corollary\label{char.surf}
Assume that $\CP$ is well defined. Then an element $u\in H_2(F)$
realizes $Du_2(\bX)$ if and only if $\CP(x)=2(u\circ x)\bmod4$ for
all $x\in\CFF0^2$.
\endcorollary

\subsection{Pontrjagin-Viro form and Rokhlin-Guillou-Marin
forms}\label{4.2}
We still assume that $X$ is an oriented closed smooth $4$-manifold,
$c$ is smooth and orientation preserving, and $F\ne\varnothing$ has
pure dimension~$2$.  Assume also that $\CP$ is well defined; due
to~\ref{char.surface} this implies that $u_2(\bX)$ is realized by a
union of components of~$F$.

\proposition\label{RGM'}
Let $F'\subset F$ be a union of components of~$F$ such that
$\CP(x)=2([F']\circ x)\bmod4$ for all $x\in\CFF0^2$.
Let $H'=H_1(F')\cap\hat\CFF1^2$ and
define a quadratic function $\CP'\:H'\to\Z/4$
via $x_1\mapsto\CP(x_1+x_0)+2([F']\circ x_0)$, where
$x_0\in H_0(F)$ is any element such that $x_1+x_0\in\CF^2$.
Then $\CP'$ coincides with the Rokhlin-Guillou-Marin form~$\gm'$ of the
characteristic surface~$F'$ in~$\bX$. In particular, $(H',\CP')$ is an
informative subspace of $H_1(F')$.
\endproposition

\proof
First notice that $q'$ is well defined and, due to~\ref{im1-2}, its
domain coincides with that of~$\gm'$. Pick an element $x\in H'$ and
consider a membrane~$\gM$ as in~\ref{2.2}. Let $\gM'=\pr^{-1}\gM$; it
is a closed $c$-invariant surface in~$X$. The index $\ind\gM$
(see~\ref{2.2}) equals the normal Euler number $\gM'\circ\gM'$.
(Indeed, $2\ind\gM$ is the obstruction to existence of a normal line
field on~$\gM'$; it is twice as big as the obstruction to existence
of a normal vector field.) The intersection points of $\sint\gM$
and~$F$ correspond to isolated fixed points of $c|_{\gM'}$, and all
the $1$-dimensional components of $\Fix c|_{\gM'}$ are two-sided
in~$\gM'$. The statement follows now from comparing the definitions
of~$\gm'$ and~$q'$ and Lemma~\ref{Pbv2}. (Note that the total number
of intersection points of $\sint\gM$ and~$F$ is even and, hence, so
is $\chi(\gM')$.)
\endproof

\theorem\label{Cong'}
If $F'$ and $\CP'$ are as in~\ref{RGM'}, then
$$
F'\circ F'+\Br\CP'\equiv\frac14\bigl[F\circ F+\sigma(X)\bigr]\bmod8.
$$
\endtheorem

\proof
The statement follows from Proposition~\ref{RGM'},
Theorem~\ref{RGMcong} applied to $F'\subset\bX$, and the well known
calculation of the ingredients of~\ref{RGMcong}:  the
self-intersection numbers of~$F'$ in~$X$ and~$\bX$ are related via
$(F'\circ F')_{\smash{\bX}}=2(F'\circ F')_X$, and the signature
of~$\bX$ is given by the Hirzebruch formula
$\sigma(X)=2\sigma(\bX)-F\circ F$.
\endproof

\theorem\label{RGM}
The restriction~$\CPP1$ of~$\CP$ to~$\CFF1^2$ coincides with the
Rokhlin-Guillou-Marin form~$\gm$ of the characteristic surface~$F$
in~$X$. In particular, $(\CFF1^2,\CPP1)$ is an informative subspace of
$H_1(F)$ and
$$
F\circ F+2\Br\CPP1\equiv\sigma(X)\bmod16.
$$
\endtheorem

\proof
The two forms are compared as in the previous proof: $\gm$ is
calculated via a generic membrane~$\gM$ as in~\ref{2.2} and~$\CPP1$,
via $\gM'=\gM\cup c(\gM)$. We may assume that $\gM'$ is an immersed
surface. It realizes $\bv_2(\partial\gM)$, as all the components of
$\partial\gM$ are two-sided in~$\gM'$ and the intersection points of
$\sint\gM\cap F$ are {\bf not} fixed points of the lift of~$c$ to the
normalization of~$\gM'$. Note also that the self-intersection points
of~$\gM'$ which are not on~$F$ appear in pairs and thus do not
contribute to $\gM'\circ\gM'$.
\endproof

\remark{Remark}
Let $X$ be the complexification of a real algebraic surface and
$c=\conj$ the Galois involution on~$X$. (More generally, one can
assume that $X$ is a compact smooth complex analytic surface and $c$
is an anti-holomorphic involution.) Then $\Fix c=\xr$ is the real
part of~$X$ and multiplication by~$\sqrt{-1}$ establishes an
isomorphism $\tau\xr=\nu\xr$. In particular, for any component~$F_i$
of~$\xr$ one has $F_i\circ F_i=-\chi(F_i)$, and the congruences
of~\ref{Cong'} and~\ref{RGM} take the form
$$
\gather
\chi(F')\equiv\frac14\bigl[\chi(\xr)-\sigma(X)\bigr]+
\Br\CP'\bmod8,\eqtag\label{RGM-alg'}\\
\chi(\xr)\equiv2\Br\CPP1-\sigma(X)\bmod16.\eqtag\label{RGM-alg}
\endgather
$$
Since $\chi(F')\equiv\Br\CP'\bmod2$, \eqref{RGM-alg'} implies
$$
\chi(\xr)\equiv\sigma(X)\bmod8.\eqtag\label{A-congr}
$$
(Certainly, \eqref{A-congr} follows as well from the Arnol$'$d
congruence for real algebraic surfaces with $[\xr]=D_Xu_2(X)$ in
$H_2(X)$.)
\endremark

\section{Real Enriques surfaces}\label{s5}

\subsection{Real Enriques surfaces}\label{5.1}
Recall that an algebraic surface~$X$ is called a \emph{$K3$-surface}
if $\pi_1(X)=0$ and $c_1(X)=0$. An algebraic surface~$E$ is called an
\emph{Enriques surface} if $\pi_1(E)=\Z/2$ and the universal
covering~$X$ of~$E$ is a $K3$-surface. (The classical definition of
Enriques surfaces is $c_1(E)\ne0$, $2c_1(E)=0$, and the relation to
$K3$-surfaces follows from the standard classification.) All
$K3$-surfaces form a single deformation family; they are all
diffeomorphic to a degree~$4$ surface in~$P^3$. Similarly, all
Enriques surfaces form a single deformation family and are all
diffeomorphic to each other.  The intersection forms of $K3$- and
Enriques surfaces are, respectively, $H_2(X;\Z)\cong3E_8\oplus2U$ and
$H_2(E;\Z)\cong E_8\oplus U$, where $E_8$ is the even unimodular form
of signature~$-8$ and $U$ is the hyperbolic plane.

A \emph{real Enriques surface} is an Enriques surface~$E$ supplied
with an anti-holo\-mor\-phic involution $\conj\:E\to E$, called
\emph{real structure}. The fixed point set $\er=\Fix\conj$ is called
the \emph{real part} of~$E$. (Obviously, these definitions apply to
any algebraic variety.)

Fix a real Enriques surface~$E$ and denote by $p\:X\to E$ its
universal covering and by $\tau\:X\to X$, the \emph{Enriques
involution} (i.e., deck translation of~$p$). The real
structure~$\conj$ on~$E$ lifts to two real structures
$\c1,\c2\:X\to X$, which commute with each other and with~$\tau$. Let
$\xri=\Fix\ci$, $i=1,2$, be their real parts. The projections
$\eri=p(\xri)$ are called \emph{halves} of~$\er$. It is easy to see
that $\er1$ and~$\er2$ are disjoint, $\er=\er1\cup\er2$, and
both~$\eri$ consist of whole components of~$\er$. Furthermore, two
components $F_1,F_2\subset\er$ belong to the same half if and only if
$\bv_1\<F_1-F_2\>=0$ (see~\cite{DK2}).

A real Enriques surface is said to be of of \emph{hyperbolic},
\emph{parabolic}, or \emph{elliptic type} if the minimal Euler
characteristic of the components of~$\er$ is negative, zero, or
positive, respectively.

\subsection{The Pontrjagin-Viro form on a real Enriques
surface}\label{5.2}
Fix a real Enriques surface~$E$. For the topological types of the
connected components of~$\er$ we will use the notation $S=S^2$,
$S_p=\#_p(S^1\times S^1)$, and $V_p=\#_p\Rp2$. The decomposition
of~$\er$ into two halves will be designated via
$\er=\{\er1\}\+\{\er2\}$.

$E$ is said to be of \emph{type~$\I{}$} if $[\er]=0$ in
$(H_2(E;\Z)/\Tors)\otimes\Z/2$ or, equivalently, $[\xr1]+[\xr2]=0$
in $H_2(X)$; otherwise $E$ is said to be of \emph{type~$\II$}.
Type~$\I{}$ is further subdivided into~$\I0$ and~$\Iu$ depending on
whether $[\er]=0$ or $Du_2(E)$ in $H_2(E)$.

\lemma\label{sufficient}
The following are sufficient conditions for the existence of the
Pontrjagin-Viro form $\CP\:\CF^2\to\Z/4$ on a real Enriques
surface~$E$\rom:
\roster
\item\local1
$E$ is an $M$-surface\rom;
\item\local2
$E$ is of type~$\Iu$ and either $\er$ is nonorientable or both~$\er1$
and~$\er2$ are nonempty\rom;
\item\local3
$E$ is of type~$\I{}$, $\er$ is nonorientable, and either both~$\er1$
and~$\er2$ are nonempty or $\er$ contains a nonorientable component
of odd genus.
\endroster
\endlemma

\proof
As shown in~\cite{DK2}, the subgroup $\CF^2\subset H_*(\er)$ is
generated by elements of the form~$[F_0]$ and $\<F_1-F_2\>\oplus
x_1$, where $F_0$, $F_1$, $F_2$ are components of~$\er$,
$x_1\in H_1(\er)$, and either $x_1^2=1$ and $F_1$, $F_2$ are in
distinct halves, or $x_1^2=0$ and $F_1$, $F_2$ are in the same half.
In particular, $E$ is Galois maximal if and only if $\er$ is
nonorientable or both~$\er1$, $\er2$ are nonempty. If $E$ satisfies
the hypotheses of~\loccit3, then $Du_2(E)\ne0$ in $\iH_2$. Since
$[\er]=Du_2(E)$ in $\iH_2$, type~$\I{}$ implies~$\Iu$. All the
statements follow now from Corollary~\ref{M-case}.
\endproof

From now on we assume that $\CP$ is defined. Two components~$F_1$,
$F_2$ of the same half are said to be in one \emph{quoter} if
$\CP\<F_1-F_2\>=0$. Since $\CP$ is linear on~$\CFF0^2$, each
half~$\eri$ splits into two quoters, which consist of whole
components of~$\eri$.  We denote this by
$\eri=(\text{quoter 1})\+(\text{quoter 2})$.
Following~\cite{Mikhalkin}, the decomposition of~$\er$ into four
quoters is called \emph{complex separation}. Due to~\ref{char.surf}
it has the following geometrical meaning: a subsurface $F'\subset\er$
is characteristic in $E/\conj$ if and only if it is the union of two
quoters which belong to distinct halves.

Let $\er=\{(\QR11)\+(\QR12)\}\+\{(\QR21)\+(\QR22)\}$ be the
decomposition of~$\er$ into quoters. If both the halves are nonempty,
denote by~$\qr1{ij}$ and~$\qr2{ji}$ the restriction to
\smash{$H_1(\QR1i)$} (respectively, \smash{$H_1(\QR2j)$}\,) of the
Rokhlin-Guillou-Marin form of the characteristic surface
$\QR1i\cup\QR2j$. As follows from~\ref{RGM'},
$$
\gather
\qr1{i1}=\qr1{i2}+Dw_1(\QR1i)\quad\text{and, hence,}\quad
\Br\qr1{i1}=-\Br\qr1{i2}\rlap,\\
\qr2{j1}=\qr2{j2}+Dw_1(\QR2i)\quad\text{and, hence,}\quad
\Br\qr2{j1}=-\Br\qr2{j2}
\endgather
$$
(see~\ref{Br} and~\ref{q-space}). If one of the halves, say, $\er2$,
is empty, denote by~$\qr1i$ the restriction of~$\CP$ to $H_2(\QR1i)$;
this form is defined on the annihilator of $w_1(\QR1i)$ and is
informative. In this notation congruence~\eqref{RGM-alg'} takes the
following form:

\proposition\label{E-RGM}
If both the halves are nonempty, then for $i,j=1,2$
$$
\chi(\QR1i)+\chi(\QR2j)\equiv
2+\tfrac14\chi(\er)+\Br\qr1{ij}+\Br\qr2{ji}\bmod8
$$
If $\er2=\varnothing$, then for $i=1,2$
$$
\chi(\QR1i)\equiv2+\tfrac14\chi(\er)+\Br\qr1i\bmod8.
$$
\endproposition

Another invariant used in the classification is the value $\CP(w_1)$,
where $w_1$ is the characteristic class of a nonorientable component
of~$\er$ of even Euler characteristic.
%It does not depend on the
%choice of a component due to the following lemma:

\proposition
Let $F_1,F_2\subset\er$ be two nonorientable components of even Euler
characteristic. Then $\CP(w_1(F_1))=\CP(w_1(F_2))$.
\endproposition

\proof
As shown in~\cite{DK1}, if $\er$ has two nonorientable components of
even Euler characteristic, it has no other nonorientable components.
Hence, $w_1(\er)=w_1(F_1)+w_1(F_2)$. On the other hand, $w_1(\er)$ is
a characteristic element in~$\CF^2$ and, due to~\iref{Br}2
and~\ref{congF}, $\CP(w_1(\er))=0$.
\endproof

\section{The Pontrjagin-Viro form via Donaldson's trick}\label{s6}

\subsection{Constructing surfaces via Donaldson's trick}\label{5.3}
Let $\tZ$ be a rational surface with real structure $c\:\tZ\to\tZ$
and nonempty real part, and $P,Q\subset\tZ$ a pair of nonsingular
real curves.

\definition[Assumption]\label{basic}
In this and next sections we assume that
\roster
\item\local1
$[P]$ and~$[Q]$ are even in $H_2(\tZ;\Z)$ and $[P]+[Q]=-2K_{\tZ}$;
\item\local2
$\dim\Ker[\operatorname{inclusion}_*\:H_2(P)\to H_2(\tZ)]=1$;
\item\local3
the multiplicity of each intersection point of~$P$ and~$Q$ is at
most~$2$.
\endroster
Since $[P]$ is even, $P_\R$ divides~$\tZ_\R$ into two parts with
common boundary~$P_\R$. Denote their closures by
$\tZ^{\pm P}=\tZ^\pm$. Similarly, introduce two parts $\tZ^{\pm Q}$
with common boundary~$Q_\R$. Let
$\tZ^{\Gd\Ge}=\tZ^{\Gd P}\cap\tZ^{\Gd Q}$ for $\Gd,\Ge=\pm$ and assume
that
\roster
\item[4]\local4
$\tZ^{++}=\varnothing$, i.e., $P_\R\subset\tZ^{-Q}$ and
$Q_\R\subset\tZ^{-P}$.
\endroster
\enddefinition

Consider the double covering $\tY\to\tZ$ branched
over~$P$ and denote by~$B$ the pull-back of~$Q$. Let $\tY'$ be the
minimal resolution of singularities of~$B$ (which occur at the
tangency points of~$P$ and~$Q$ and are all nondegenerate double
points) and~$B'$ the proper transform of~$B$. The real structure~$c$
on~$\tZ$ lifts to two real structures~$c^\pm$ on~$\tY'$. In respect
to one of them,~$c^+$, the real part~$\tY'_\R$ projects to~$\tZ^+$
and $B'_\R=\varnothing$.

\proposition\label{rational}
$\tY$ and~$\tY'$ are rational surfaces, $[B]=-2K_{\tY}$, and
$[B']=-2K_{\tY'}$. The double covering $\tX\to\tY'$ of~$\tY'$
branched over~$B'$ is a $K3$-surface.
\endproposition

\proof
The relations $[B]=-2K_{\tY}$ and $[B']=-2K_{\tY'}$ follow from the
projection formula and~\iref{basic}1. Thus, the anti-bicanonical
class of~$\tY$ is effective and, hence, $\tY$ is either rational or
ruled. On the other hand, from the Smith exact sequence and~\ditto2
it follows that $H_1(\tY)=0$. Hence, $\tY$ is rational. Now the
projection formula gives $2K_{\tX}=0$, i.e., $\tX$ is
a minimal surface of Kodaira dimension~$0$. Using the Riemann-Hurwitz
and adjunction formulas one obtains
$\chi(\tX)=2\chi(\tY')+2K^2_{\tY'}=24$. Hence, $\tX$ is a
$K3$-surface.
\endproof

Denote by~$\c1$ the deck translation of $\tX\to\tY'$. Due to~\loccit4
above one of the two lifts of~$c^+$ to~$\tX$ is fixed point free;
denote it by~$\tau$. One can now apply to $(\c1,\tau)$ the following
equivariant version of Donaldson's trick:

\proposition[Proposition \rm(cf.~\cite{DK3})]\label{Donaldson}
Let~$\tX$ be a $K3$-surface and $(c_h,c_a)$ a pair of commuting
involutions on~$\tX$, one holomorphic and one anti-holomorphic. Then
there is a complex structure on~$\tX$ in respect to which $c_h$ is
anti-holomorphic and $c_a$ is holomorphic.
\endproposition

Let~$\ttX$ be the resulting $K3$-surface. Then the quotient
$E=\ttX/\tau$ is an Enriques surface and $\c1$ descends to a real
structure~$\conj$ on~$E$. Clearly, $\er1=B'/c^+$ and $\er2=\tY'_\R$,
and there is a projection $\pi\:\er\to\bQ\cup\tZ^+$, which is a
branched double covering outside the tangency points of~$P_\R$
and~$Q_\R$.  The pull-back~$\pi^{-1}(T)$ of each tangency point~$T$
consists of a one-sided loop in~$\er2$ and a point in~$\er1$.
Let $\rho\:\er\to\er$ be the deck translation of~$\pi$.

\lemma[Notation]\label{topER1}
For a subset $S\subset\bQ\cup\tZ^+$ we denote by
$\<\pi^{-1}(Q)\>\in H_0(\er)\>$ the class generated by the connected
components of~$\pi^{-1}(S)$, and by $[\pi^{-1}(S)]\in H_*(\er)$, the
fundamental class of its components of highest dimension
\rom(provided that they are all closed manifolds\rom).
\endlemma

\subsection{Existence of~$\CP$ and complex separation}\label{5.4}
Let $E$ be a real Enriques surface obtained as in~\ref{5.3} from a
configuration $(\tZ;P,Q)$. For all intermediate objects we keep
the notation introduced in~\ref{5.3}.

A nonsingular real curve~$C\subset\tZ$ is said to be of
\emph{type~\rom{I}}, or \emph{separating}, if $C/c$ is orientable.
The real part of each real component~$C_i$ of~$C$ has a
distinguished pair of opposite orientations, called \emph{complex
orientations}, which are induced from an orientation of~$C_i/c$.
The complement $C_i\sminus C_{i,\R}$ consists of
two components~$C_i^\pm$ with common boundary $C_{i,\R}$; their
natural orientations induce the two complex orientations of
$C_{i,\R}$.

A triple $(\tZ;P,Q)$ as in~\ref{5.3} is said to be of
\emph{type~\rom{I}} if $P$ is of type~I and $\tZ^-\sminus Q_\R$
is orientable. Let $P_i$, $i=1,\dots,p$, and $Q_j$, $j=1,\dots,q$,
be the real components of the curves, $Z^-_k$, $k=1,\dots,z^-$, the
connected components of $\tZ^{-+}$ and $\tZ^{--}$, and $Z^+_l$,
$l=1,\dots,z^+$, the connected components of $\tZ^+$. Fix some
orientations of~$Z^-_k$ and some complex orientations of~$P_{i,\R}$;
this determines an orientation of $\bP_i$ and a distinguished
half~$P^+_i$ for each real component~$P_i$. Thus, the fundamental
classes $[\bP_i]$ and $[Z^-_k]$ are well defined $\!{}\bmod4$ in all
homology groups where they make sense. The classes $[\bQ_j]$ and
$[Z^+_l]$ are defined $\!{}\bmod2$.

\definition
We say that a triple $(\tZ;P,Q)$ of type~\rom{I} \emph{admits a
fundamental cycle} if there are some odd integers~$\Gl_i$,
$i=1,\dots,p$, and $\Gk^-_k$, $k=1,\dots,z^-$, and some
integers~$\Gm_j$, $j=1,\dots,q$, and~$\Gk^+_l$, $l=1,\dots,z^+$ such
that
$$
\tsize
\sum\Gl_i[P_{i,\R}]+\sum\Gk^-_k[\partial Z^-_k]=
2\sum\Gm_j[Q_{j,\R}]+2\sum\Gk^+_l[\partial Z^+_l]
\eqtag\label{fund.def}
$$
in $H_1(P_\R\cup Q_\R;\Z/4)$. The combination
$$
\tsize
\gC=\sum\Gl_i[\bP_i]+\sum\Gk^-_k[Z^-_k]+
2\sum\Gm_j[\bQ_j]+2\sum\Gk^+_l[Z^+_l],
\eqtag\label{cycle}
$$
which is a $(\!\!{}\bmod4)$-cycle in~$\bZ$, is called a
\emph{fundamental cycle}. It is called
\emph{proper} if $[\gC]=2Dw_2(\bZ)$ in $H_2(\bZ;\Z/4)$.
\enddefinition

\proposition\label{good.chain}
A triple $(\tZ;P,Q)$ of type~\rom{I} admits a fundamental cycle if and
only if $\sum[P_{i,\R}]$ belongs to the subgroup in $H_1(\tZ^-;\Z/4)$
spanned by the classes $2[P_{i,\R}]$, $2[Q_{j,\R}]$, and
$2\partial[Z^+_l]$. If $\tZ$ is an $M$-surface with $\tZ_\R$
connected, any fundamental cycle is proper.
\endproposition

\proof
The first part follows from the exact sequence of pair
$(\tZ^-,P_\R\cup Q_\R)$. If $\tZ$ is an $M$-surface with $\tZ_\R$
connected, then $\bZ$ is a $\Z$-homology sphere (see~\ref{sphere})
and $[\gC]=2w_2(\bZ)$ holds trivially.
\endproof

%Assume that the resulting real Enriques surface~$E$ has well-defined
%Pontrjagin-Viro form.  Its complex separation is called
%\emph{$\rho$-invariant} if all quoters are fixed by~$\rho$.

The Pontrjagin-Viro form on a real Enriques surface~$E$ is said to
have \emph{$\rho$-invariant complex separation} if all quoters are
fixed by~$\rho$. Since $\CP$ is obviously $\rho_*$-invariant, for
each half~$\eri$, $i=1,2$, one has either $\rho(\QR{i}j)=\QR{i}j$,
$j=1,2$, or $\rho(\QR{i}1)=\QR{i}2$. This remark is sufficient to
eliminate the possibility of noninvariant complex separation in all
cases considered in Section~\ref{s7} below.

\theorem\label{main.th}
The real Enriques surface resulting from a triple $(\tZ;P,Q)$ has
Pontrjagin-Viro form with $\rho$-invariant complex separation if and
only if all tangency points of~$P$ and~$Q$ are real and $(\tZ;P,Q)$
is of type~\rom{I} and admits a proper fundamental cycle. If this is
the case, the complex separation is determined by a proper
fundamental cycle~\eqref{cycle}\rom:
\roster
\item
the components $\pi^{-1}(\bQ_a)$, $\pi^{-1}(\bQ_b)$ corresponding to
real components $Q_a,Q_b$ of~$Q$ belong to the same quoter if and
only if $\Gm_a-\Gm_b\equiv0\bmod2$\rom;
\item
the components~$\pi^{-1}(Z^+_a)$, $\pi^{-1}(Z^+_b)$ corresponding to
$Z^+_a,Z^+_b\subset\tZ^+$ belong to the same quoter if
and only if \smash{$\Gk^+_a-\Gk^+_b\equiv0\bmod2$}.
\endroster
\endtheorem

\corollary\label{unique.cycle}
A proper fundamental cycle~\eqref{cycle} is unique $\!{}\bmod4$ up to
$2[\tZ_\R]$, $2([\tZ^+]+[\bP])$, and $2([\tZ^{-+}]+[\bQ])$.
\endcorollary

\subsection{Proof of Theorem~\ref{main.th}}\label{5.5}
Note that the construction of~\ref{5.3} still works if \iref{basic}1
and~\ditto2 are replaced with a weaker assumption
\widestnumber\item{$(1')$}
\roster
\item"$(1')$"\label{weak}
$[P]$ and $[Q]$ are divisible by~$2$ in $H_2(\tZ;\Z)$.
\endroster
Certainly, Proposition~\ref{rational} does not hold in this case and
\ref{Donaldson} does not apply; thus, $E$ is just a $4$-manifold with
orientation preserving involution. In view of~\ref{char.surface}
Theorem~\ref{main.th} would follow from the following more general
result:

\theorem\label{general.main}
Let $(E,\conj)$ be the manifold with involution resulting from a
triple $(\tZ;P,Q)$ satisfying~\ref{weak} above and~\iref{basic}3,
\ditto4. Then $\er\subset E/\conj$ contains a $\rho$-invariant
characteristic surface if and only if all tangency points of~$P$
and~$Q$ are real and $(\tZ;P,Q)$ is of type~\rom{I} and admits a
proper fundamental cycle. If this is the case, the $\rho$-invariant
characteristic subsurfaces of~$\er$ are those of the form
$\sum\Gm_j[\pi^{-1}(\bQ_j)]+\sum\Gk^+_l[\pi^{-1}(Z^+_l)]$, where
$\Gm_j$ and~$\Gk^+_l$ are the coefficients of a proper fundamental
cycle.
\endtheorem

In order to prove Theorem~\ref{general.main} we replace $E/\conj$
with $\bY'=\tY'/c^+$. In the construction $\tY'$ is obtained
from~$\tY$ by blow-ups of the tangency points of~$P$ and~$Q$. If such
a point is real, the blow-up results in connected summation of~$\bY$
and $\smash{\overline{\Cp2}/\text{conjugation}}\cong S^4$ and thus
does not affect the topology. A pair of conjugate tangency points
results in the (topological) blow-up of their common image in~$\bY$.
This produces a $(-1)$-sphere in~$\bY'$ which is
$(\!\!{}\bmod2)$-orthogonal to all pull-backs $[\pi^{-1}(\bP_j)]$ and
$[\pi^{-1}(Z^+_l)]$, which shows that $\er$ does not contain a
$\rho$-invariant characteristic surface. Thus, we can assume that $P$
and~$Q$ do not have imaginary tangency points and replace~$\bY'$
with~$\bY$.

Clearly, $\bY$ is the double covering of~$\bZ=\tZ/c$ branched over the
\emph{Arnol$\,'\!$d surface} $\gA^-=\bP\cup\tZ^-$. Denote by
$\pr\:\bY\to\bZ$ the projection and by
$\omega\in H^1(\bZ\sminus\gA^-)$ its characteristic class. Since
$\partial\:H_3(\bZ,\gA^-)\to H_2(\gA^-)$ is a monomorphism, $\omega$
is uniquely characterized by the property $\partial D\omega=[\gA^-]$.
Note also that, since $\bZ$ is orientable,
$\omega\scap D\omega=\Sq_1D\omega\in H_2(\bZ,\gA^-)$ and
$\partial\Sq_1D\omega=\Sq_1[\gA^-]=Dw_1(\gA^-)$.

\lemma\label{6.3.2}
$\sum\Gm_j[\pi^{-1}(\bQ_j)]+\sum\Gk^+_l[\pi^{-1}(Z^+_l)]=Dw_2(\bY)$
in $H_2(\bY)$ if and only if in $H_2(\bZ,\gA^-)$
$$
\tsize
\sum\Gm_j[\bQ_j]+\sum\Gk^+_l[Z^+_l]=\Sq_1D\omega+Dw_2(\bZ).
\eqtag\label{rel.proper}
$$
\endlemma

\proof
The statement follows immediately from the projection formula and
Smith exact sequence.
\endproof

\lemma\label{6.3.4}
A linear combination~\eqref{cycle} is a fundamental cycle if and only
if $\sum\Gm_j[\partial\bQ_j]+\sum\Gk^+_l[\partial Z^+_l]=Dw_1(\gA^-)$
in $H_1(\gA^-)$.
\endlemma

\proof
This is a direct consequence of the definition of Bockstein
homomorphism via chains.
\endproof

From Lemma~\ref{6.3.4} it follows that a necessary condition for
\eqref{rel.proper} to hold is that $\Gm_j$, $\Gk^+_l$ must be
coefficients of a fundamental cycle; in particular, this implies that
$(\tZ;P,Q)$ must be of type~I. If this is the case, $\Sq_1D\omega$ is
the relativization of a class $x\in H_2(\bZ,\bP_\R\cup\bQ_\R)$,
which is well defined up to the image of $H_2(\gA^-)$ and has the
property $2x=\sum[\bP_i]+\sum[Z^-_k]\bmod2H_2(\gA^-)$ in
$H_2(\bZ,\bP_\R\cup\bQ_\R;\Z/4)$. Since both
$\stimes2\:H_2(\bZ,\bP_\R\cup\bQ_\R)\to
H_2(\bZ,\bP_\R\cup\bQ_\R;\Z/4)$
and relativization
$H_2(\bZ;\Z/4)\to H_2(\bZ,\bP_\R\cup\bQ_\R;\Z/4)$ are
monomorphisms, \eqref{rel.proper} is equivalent to the fact that
$\Gm_j$ and~$\Gk^+_l$ are coefficients of a proper fundamental cycle.
\qed

\subsection{Values on $1$-dimensional classes}\label{5.6}
Fix a triple $(\tZ;P,Q)$ satisfying the conditions of
Theorem~\ref{main.th}, so that the Pontrjagin-Viro form $\CP$ on~$E$
is well defined and the complex separation is $\rho$-invariant.

Denote $\gS=\tZ_\R\cup\bP\cup\bQ$ and for an immersed (in the
obvious sense) loop $\gtl\subset\gS$ transversal to~$P_\R$ and~$Q_\R$
define its `normal Euler number' $e(\gtl)\in\Z/2$ to be~$1$ or
$0\bmod2$ depending on whether $\gtl$ is disorienting or not. (If
$\gtl$ passes through an isolated intersection point of $\bP$ and
$\bQ$, the orientation is transferred so that the point have
intersection index~$+1$.) The following obvious observation is
helpful in evaluating $e(\gtl)$: {\sl let an oriented arc~$\gtl'$
belong to a half~$C^+$ of a type~\rom{I} curve~$C$ \rom(which, in
our case, can be union of real components of~$P$ and separating real
components of~$Q$\rom) so that $\partial\gtl'\subset C_\R$ and
$\gtl'$ is normal to~$C_\R$. Then the co-orientation induced from the
complex orientation of~$C_\R$ at the initial point of~$\gtl'$ is
transferred by~$\gtl'$ to the co-orientation opposite to the complex
orientation of~$C_\R$ at the terminal point of~$\gtl'$.}

\proposition\label{pull-back}
For $\gtl\subset\gS$ as above one has
$$
\tsize
\CP([\pi^{-1}(\gtl)]\oplus\<\pi^{-1}(\gtl)\>)=
2e(\gtl)+2i_Q(\gtl)+2i^+(\gtl)+i_{P\cap Q}(\gtl)\bmod4,
$$
where $i_Q(\gtl)$ and $i^+(\gtl)$ are the numbers of isolated
intersection points of~$\gtl$ with~$\bQ$ and $\tZ^+$, respectively,
and  $i_{P\cap Q}(\gtl)$ is the number of intersection points
$\bP\cap\bQ$ through which $\gtl$ passes.
\endproposition

\proof
Let $\gM$ be a membrane in~$\bZ$ normal along~$\partial\gM=\gtl$
and transversal in $\sint\gM$ to all strata of~$\gS$. Then
the value in question can be found using
$\gM_Y=\pr^{-1}(\gM)$ (cf.~\ref{RGM'}). Clearly,
$\ind\gM_Y=2\ind\gM+\frac12i_{P\cap Q}(\gtl)$. (To define
$\ind\gM$, we use a normal line field on~$\gtl$ tangent to all
strata of~$\gS$ and patch it at the points of $\bP\cap\bQ$ using
local orientations.) Further,
$\Card(\sint\gM_Y\cap\er1)=i_Q(\gtl)+2\Card(\sint\gM\cap\bQ)$ and
$\Card(\sint\gM_Y\cap\er2)=i^+(\gtl)+2\Card(\sint\gM\cap\tZ^+)$.  It
remains to notice that $2\ind\gM=0$ or $1\bmod2$ depending on whether
the line field is orientable or not.
\endproof

\proposition\label{tang.pt}
For a point $T\in P_\R\cap Q_\R$ one has
$\CP([\pi^{-1}(T)]+\<\pi^{-1}(T)\>)=1$.
\endproposition

Let~$S$ be a connected component of one of $\bQ$, $\bP$, or
$\tZ_\R\sminus(P_\R\cup Q_\R)$. Assume that $\partial S$ does not
contain a tangency point of~$P$ and~$Q$. Then $\partial S$ is a loop
in~$\er$.

\proposition\label{Z-}
Let $Z^-_k$ be a connected component of~$\tZ^-$ with $\partial Z^-_k$
disjoint from $P_\R\cap Q_\R$. Then
$\CP[\partial Z^-_k]=2\chi(Z^-_k)\bmod4$.
\endproposition

\proposition\label{P}
Let $P_i$ be a real component of~$P$ with $P_{i,\R}$ disjoint from
$P_\R\cap Q_\R$.
Then
$\CP([P_{i,\R}]+\<\pi^{-1}(\bP_i\cap\bQ)\>)=\frac12[P_i]^2\bmod4$.
\endproposition

\proof[Proof of Propositions~\ref{tang.pt},~\ref{Z-},
and~\ref{tang.pt}]
We lift~$Z^-_k$ (respectively,~$P_i$ or the
exceptional curve appearing when~$T$ is blown up in~$\tY$) to~$\tX$
and then project the result to~$E$.  This gives a $\conj$-invariant
closed surface in~$E$, and the value in question is found
via~\ref{Pbv2}.
\endproof

In the rest of this section we consider the case when $S$ is
either~$\bQ_j$ for a real component~$Q_j$ of~$Q$ or a
component~$Z^+_l$ of~$\tZ^+$. Fix a proper fundamental cycle~$\gC$
and denote by~$U$ the corresponding characteristic surface in~$\bY$.
Assume that $U$ contains $\pi^{-1}(S)\subset\er$ and denote by
$\gm_{\gC}$ the Rokhlin-Guillou-Marin form of~$U$. The choice
of~$\gC$ determines a preferred orientation of~$\bP$ (which induces
$\sum\Gl_i[P_{i,\R}]\bmod4$ on $\partial\bP$) or, equivalently, a
half~$P^+$ of $P\sminus P_\R$.

Let $\omega\in H^1(\tZ_\R)$ be the class Poincar\'e dual to~$[C_\R]$,
where $C\subset\tZ$ is a real curve with $2[C]=[P]$ in $H_2(\tZ;\Z)$.
The restriction of~$\omega$ to $\tZ_\R\sminus P_\R$ is the
characteristic class of the restricted covering $\tY\to\tZ$
(see~\ref{char.class}). For each real component~$Q_j$ of~$Q$
denote by $\omega_j\in H^1(\bQ_j\sminus\bP)$ the characteristic class
of the covering $\bY\to\bZ$ restricted to $\bQ_j\sminus\bP$.
($\omega_j$ can be interpreted as the linking number with~$\gA^-$
in~$\bZ$.)

{\it Assume that $\<\omega,[\gtl]\>=0$ for each boundary component
$\gtl\subset\partial S$} and define the `linking number'
$\lk_{\gC}w_1(S)\in\Z/4$ of the characteristic class of~$S$ with the
Arnold surface. Fix some orientations of the boundary components and,
if $S=\bQ_j$, some local orientations at the intersection
points $S\cap\bP$. This defines a lift of~$w_1(S)$ to a class
$w_1'\in H^1(S,\partial S\cup(S\cap\bP))$ (which can be defined as
the obstruction to extending the chosen orientations to the whole
surface). We let $\lk_{\gC}w_1(S)=2\<\omega,Dw_1'\>$ if $S=Z^+_l$ or
$2\<\omega_j,Dw_1'\>-(\sint S\circ\bP)$ if $S=\bQ_j$. (In the latter
case the intersection index is defined~$\!{}\bmod4$ using the chosen
local orientations of~$S$; the condition
$\<\omega_j,[\partial S]\>=0$ implies that $\partial S$ is not linked
with~$\gA^-$ and, hence, $\sint S\circ\bP=0\bmod2$.)

In the case $S=\bQ_j$ for a real component~$Q_j$ of~$Q$ the above
definition is cumbersome and not `seen' in the real part. If $Q_j$ is
of type~I, it can be simplified:
$\lk_{\gC}w_1(S)=\Card(Q_j^+\cap P^+)-\Card(Q_j^+\cap P^-)\bmod4$ for
any half~$Q_j^+$ of~$Q_j$.

\proposition\label{Q,Z+}
Let $S$ be either~$\bQ_j$ for a real component~$Q_j$ of~$Q$,
or a component~$Z^+_l$ of~$\tZ^+$. Assume that $\partial S$ is
disjoint from $P_\R\cap Q_\R$ and $\<\omega,[\gtl]\>=0$ for each
boundary component $\gtl\subset\partial S$. Then
$\gm_{\gC}[\partial S]=\lk_{\gC}w_1(S)$.
\endproposition

\remark{Remark}
Proposition~\ref{Q,Z+} applies to the generalized construction
described in~\ref{5.5}. In the case of Enriques surfaces, due to
Corollary~\ref{unique.cycle}, the preferred orientation of~$\bP$ is
defined up to total reversing and $\lk_{\gC}w_1(S)$ does not depend
on the choice of~$\gC$. Furthermore, due to the assumption made
$\partial S$ is a collection of two-sided circles in~$\er$. Hence,
$\gm_{\gC}[\partial S]=\CP[\partial S]$.
\endremark

\proof
Let $\gM_Y\subset\bY$ be an oriented membrane normal along
$\partial\gM_Y=\partial S$ and transversal in $\sint\gM$ to both~$U$
and~$\gA^-$.  Let $\gM$ be its projection to~$\bZ$.  Clearly,
$\ind\gM_Y=\ind\gM$ and the intersection points of $\sint\gM_Y$
and~$U$ project one-to-one.  Furthermore, $\gM$ is tangent to~$\gA^-$
at its inner points; hence,
$$
2\Card(\sint\gM_Y\cap U)=
(\sint\gM\circ\gC)-2(\sint\gM\circ\gA^-)\bmod4.
$$
(Recall that~$U$ consists precisely of those components of~$\er$
whose coefficients in~$\gC$ are $2\bmod4$.) Since
$[\gC]=2w_2(\bZ)\bmod4$ and $[\gA^-]=0\bmod2$ in~$\bZ$, the
expressions $2\ind\gM+(\sint\gM\circ\gC)\bmod4$ and
$(\sint\gM\circ\gA^-)\bmod2$ do not depend on the choice of~$\gM$; one
can replace~$\gM$ with another membrane, which does not have to lift
to~$\bY$. (Strictly speaking, to claim this one should fix the
orientation of~$\partial\gM$ induced from~$\gM$; however, it is
chosen arbitrarily for $\gM_Y$ and does not affect the result.)

Take for a new membrane $S$ shifted along a normal vector field. To
make it orientable, fix some choices used to define~$\lk_{\gC}w_1(S)$,
cut~$S$ along a simple loop~$\gtl$ representing~$w_1'$, pick a
generic orientable membrane~$\gN$ spanned by~$\gtl$, and attach
$2\gN$ to the cut. Let~$\gM$ be the result. Then
$\sint\gM\circ\gA^-=\Card(S\cap\bP)\bmod2$ and
$$
\sint\gM\circ\gC=
2\ind S+(\sint S\circ\gA^-)+2\Card(\gN\cap\gA^-)\bmod4.
$$
(The first term here is due to the original shift of~$S$ along a
normal field; recall that $S$ has coefficient $2\bmod4$ in~$\gC$. The
second term stands for the intersection index of~$\bP$ and the cut
of~$S$, which are both oriented.) Since $\ind\gM=\ind S\bmod2$,
one obtains
$\gm_{\gC}(\partial S)=2\Card(\gN\cap\gA^-)-(\sint S\circ\bP)=
\lk_{\gC}w_1(S)\bmod4$.
\endproof

\subsection{Resolving singularities of~$P$ and~$Q$}\label{5.7}
Let $(\tZ';P',Q')$ be a triple satisfying \iref{basic}1--\ditto4,
except that $P'$ and~$Q'$ may be singular. Assume that the curve
$P'+Q'$ has at most simple singular points (i.e., those of
type~$A_p$, $D_q$, $E_6$, $E_7$, or~$E_8$). Then there is a sequence
of blow-ups which converts $(\tZ';P',Q')$ to a triple $(\tZ;P,Q)$
with $P$ and~$Q$ nonsingular and satisfying \iref{basic}1--\ditto4.
More precisely, the singularities of $P'+Q'$ can be resolved by
blowing up double or triple points. Let $O$ be a singular point of,
say,~$P'$. Blow it up and denote by~$e$ the exceptional divisor and
by~$\tilde P$, the proper transform of~$P'$. Then the new pair
$(P,Q)$ on the resulting surface is constructed as follows: $Q$ is
the full transform of~$Q'$ and $P$ is either~$\tilde P$ or
$\tilde P+e$, depending on whether $O$ is a double or triple point
of~$P'$. The singular points of~$Q'$ are resolved similarly, with
$P'$ and~$Q'$ interchanged. If the resulting curves $(P,Q)$ are still
singular, the procedure is repeated.

%\subsection{Noninvariant separations}\label{5.8}
%Necessary conditions for existence.\mnote{to write!!!}

\section{Calculation for real Enriques surfaces}\label{s7}

\subsection{$M$-surfaces of elliptic and parabolic type}\label{7.1}

\theorem[Theorem \rm(see \cite{DIK})]\label{PV-ep}
A real Enriques $M$-surface~$E$ of parabolic or elliptic type is
determined up to deformation equivalence by its complex separation
and the value $\CP(w_1)$ of~$\CP$ on the characteristic element of a
nonorientable component of~$\er$ of even Euler characteristic \rom(if
such a component exists\rom). The deformation types of such surfaces
are given in Tables~\ref{4V+2S} and~\ref{parabolic}, which list the
separations of the two halves and the possible values of~$\CP(w_1)$.
\endtheorem

\remark{Remark}
In~\cite{Kucuk} it is shown that in all cases listed in the tables
$\CP$ is uniquely recovered (up to autohomeomorphism of~$\er$
preserving the complex separation) from the complex separation and
$\CP(w_1)$ via~\ref{E-RGM} and, moreover, all forms satisfying the
congruences of~\ref{E-RGM} are realized by real Enriques surfaces.
\endremark

\midinsert
\table\label{4V+2S}
$M$-surfaces of elliptic type ($\er=4V_1\+2S$)
\endtable
\hbox to\hsize{\eightpoint\let\0\varnothing
\hss\vtop{\halign{$#$\hss&$#$\hss&\quad$#$\hss\cr
&(2V_1\+S)\+(2V_1\+S)&\0\cr
&(4V_1)\+(2S)&\0\cr
\noalign{\smallskip}
&(2V_1\+S)\+(V_1\+S)&(V_1)\cr
&(3V_1)\+(2S)&(V_1)\cr
\noalign{\smallskip}
&(2V_1\+S)\+(S)&(V_1)\+(V_1)\cr
&(V_1\+S)\+(V_1\+S)&(V_1)\+(V_1)\cr
&(2V_1)\+(2S)&(V_1)\+(V_1)\cr
&(V_1\+S)\+(V_1\+S)&(2V_1)\cr
&(2V_1)\+(2S)&(2V_1)\cr
&(3V_1)\+(V_1\+S)&(S)\cr}}\hss\vrule
\hss\vtop{\halign{$#$\hss&$#$\hss&\quad$#$\hss\cr
&(V_1\+S)\+(S)&(2V_1)\+(V_1)\cr
&(V_1)\+(2S)&(2V_1)\+(V_1)\cr
&(3V_1)\+(S)&(V_1)\+(S)\cr
&(2V_1)\+(V_1\+S)&(V_1)\+(S)\cr
&(2V_1)\+(V_1\+S)&(V_1\+S)\cr
\noalign{\smallskip}
&(S)\+(S)&(2V_1)\+(2V_1)\cr
&(2S)&(2V_1)\+(2V_1)\cr
&(2V_1)\+(S)&(V_1\+S)\+(V_1)\cr
&(V_1\+S)\+(V_1)&(V_1\+S)\+(V_1)\cr
&(2V_1)\+(S)&(2V_1)\+(S)\cr}}\hss}
\endinsert

\midinsert
\table\label{parabolic}
$M$-surfaces of parabolic type
\endtable
\vglue-\medskipamount
\let\style\textstyle
\def\gap{\quad}
\def\\{\noalign{\smallskip}}
\catcode`\+\active \let+\+
\let\0\varnothing
\def\c{0}
\def\n{2}
\def\cn{\hidewidth0,2\hidewidth}

\hbox to\hsize{\eightpoint\hss\vtop{%
\halign{\llap{$\style^#{}\,$}&$\style#$\hss&
  \gap$\style#$\hss&\gap\hss$\style#$\hss\cr
\multispan3\hss\bf Case $E_{\R}=S_1+V_2+4S$\hss\cr
\noalign{\smallskip}
 &(V_2+2S)+(2S)&(S_1)+(\0)&\c\cr
\noalign{\medskip}
\multispan3\hss\bf Case $E_{\R}=2V_2+4S$\hss\cr
\noalign{\smallskip}
 &(V_2)+(V_2)&(2S)+(2S)&\c\cr
 &(V_2)+(V_2)&(3S)+(S)&\n\cr
 &(2V_2)+(\0)&(2S)+(2S)&\cn\cr\\
 &(V_2+2S)+(2S)&(V_2)+(\0)&\c\cr\\
 &(V_2+S)+(2S)&(V_2+S)+(\0)&\n\cr
 &(V_2+2S)+(S)&(V_2)+(S)&\n\cr
 &(V_2+S)+(S)&(V_2+S)+(S)&\c\cr
\noalign{\medskip}
\multispan3\hss\bf Case $E_{\R}=V_2+2V_1+3S$\hss\cr
\noalign{\smallskip}
 &(V_2+2S)+(2V_1+S)&\0&\c\cr
 &(V_2+2V_1+S)+(2S)&\0&\cn\cr\\
 &(V_2+V_1+S)+(V_1+S)&(S)+(\0)&\cn\cr\\
 &(V_2+S)+(2V_1)&(S)+(S)&\c\cr
 &(V_2+S)+(2V_1)&(2S)+(\0)&\n\cr
 &(V_2+2V_1)+(S)&(S)+(S)&\cn\cr\\
 &(V_2+V_1)+(V_1)&(2S)+(S)&\cn\cr\\
 &(V_2+2S)+(V_1+S)&(V_1)+(\0)&\c\cr
 &(V_2+V_1+S)+(2S)&(V_1)+(\0)&\cn\cr}}
\hss\vrule\hss\vtop{%
\halign{\llap{$\style^#{}\,$}&$\style#$\hss&
  \gap$\style#$\hss&\gap\hss$\style#$\hss\cr
\multispan4\hss\bf Case $E_{\R}=V_2+2V_1+3S$ \rm(continued)\hss\cr
\noalign{\smallskip}
 &(V_2+S)+(V_1+S)&(V_1)+(S)&\c\cr
 &(V_2+S)+(V_1+S)&(V_1+S)+(\0)&\n\cr
 &(V_2+V_1+S)+(S)&(V_1)+(S)&\cn\cr\\
 &(V_2+S)+(V_1)&(V_1+S)+(S)&\c\cr
 &(V_2+S)+(V_1)&(V_1)+(2S)&\n\cr
 &(V_2+V_1)+(S)&(V_1+S)+(S)&\cn\cr\\
 &(V_2)+(V_1)&(V_1+S)+(2S)&\c\cr
 &(V_2)+(V_1)&(V_1+2S)+(S)&\n\cr
 &(V_2+V_1)+(\0)&(V_1+S)+(2S)&\cn\cr\\
 &(V_2+S)+(2S)&(V_1)+(V_1)&\c\cr
 &(V_2+S)+(2S)&(2V_1)+(\0)&\n\cr
 &(V_2+2S)+(S)&(V_1)+(V_1)&\c\cr\\
 &(V_2+S)+(S)&(2V_1)+(S)&\c\cr
 &(V_2+S)+(S)&(V_1+S)+(V_1)&\n\cr\\
 &(V_2)+(S)&(V_1+S)+(V_1+S)&\c\cr
 &(V_2)+(S)&(2V_1+S)+(S)&\n\cr
 &(V_2+S)+(\0)&(V_1+S)+(V_1+S)&\c\cr
 &(V_2+S)+(\0)&(2V_1)+(2S)&\n\cr\\
 &(V_2)+(\0)&(2V_1+S)+(2S)&\c\cr
 &(V_2)+(\0)&(V_1+2S)+(V_1+S)&\n\cr
}}\hss}
\endinsert

The calculation of Pontrjagin-Viro forms is based on the results of
Section~\ref{s5} and the following statement, which gives explicit
models of $M$-surfaces of elliptic and parabolic types:

\theorem[Theorem \rm(see \cite{DIK})]\label{constr-ep}
Up to deformation any real Enriques $M$-surface of
parabolic or elliptic type can be obtained by the construction
of~\ref{5.3} and~\ref{5.7} from a triple $(\tZ';P',Q')$, where either
\roster
\item
$\tZ'=\Rp2$, $P'$ is an $M$-curve of degree~$4$,
and $Q'$ is a pair of lines, or
\item
$\tZ'$ is a hyperboloid $\Rp1\times\Rp1$, $P'$ is a nonsingular
$M$-curve of bi-degree $(4,2)$, and~$Q'$ is a pair of generatrices of
bi-degree $(0,1)$.
\endroster
\endtheorem

Figure~\ref{fig1} illustrates the construction of the four
nonequivalent surfaces with $\er=\{2V_2\}\+\{4S\}$ (see
Table~\ref{parabolic}). To emphasize the difference the
linear components of~$Q'$ (the lines) are shown tangent to~$P'$; in
reality they must be shifted away from~$P'$.

\midinsert
\bgroup
\eightpoint

\let\0\varnothing
\def\ctwo#1#2{\hbox to\hsize{\quad\hss#1\hss\qquad\hss#2\hss\quad}}
\ctwo{\epsfbox{hpb1.eps}}{\epsfbox{hpb2.eps}}
\kern3pt
\ctwo{\strut$\{(2V_2)\+(\0)\}\+\{(2S)\+(2S)\}$, $\CP(w_1)=0$}
  {$\{(2V_2)\+(\0)\}\+\{(2S)\+(2S)\}$, $\CP(w_1)=2$}
\hbox to\hsize{\hss\strut
  ($\tZ'$ is a hyperboloid $\Rp1\times\Rp1$)\hss}
\medskip
\ctwo{\epsfbox{plane2.eps}}{\epsfbox{plane1.eps}}
\kern3pt
\ctwo{\strut$\{(V_2)\+(V_2)\}\+\{(2S)\+(2S)\}$, $\CP(w_1)=0$}
  {$\{(V_2)\+(V_2)\}\+\{(3S)\+(S)\}$, $\CP(w_1)=2$}
\hbox to\hsize{\hss\strut
  ($\tZ'$ is a real projective plane $\Rp2$)\hss}
\egroup
\figure\label{fig1}
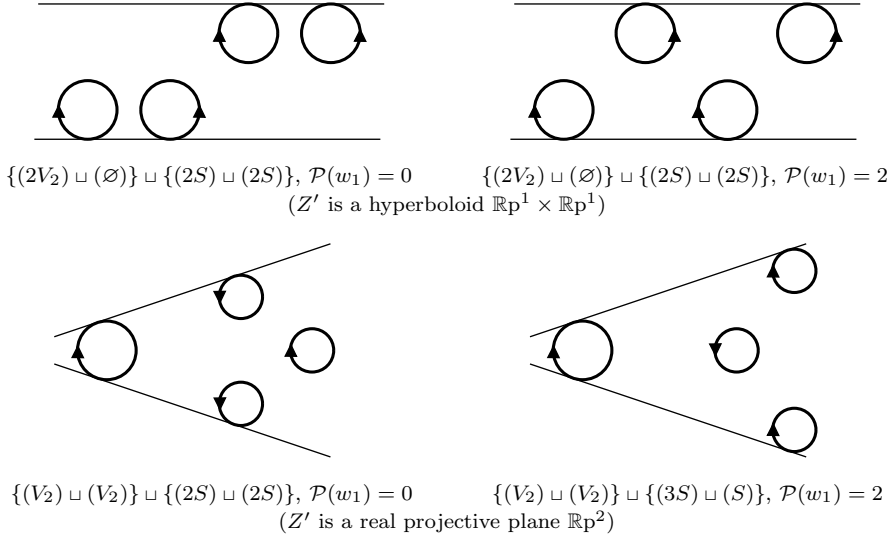
Models of real Enriques surfaces with $\er=\{2V_2\}\+\{4S\}$
\endfigure
\endinsert

\subsection{$M$-surfaces of hyperbolic type}\label{7.2}

\theorem[Theorem \rm(see \cite{DK3})]\label{PV-hyp}
A real Enriques surface of hyperbolic type is determined up to
deformation by the decomposition $\er=\{\er1\}\+\{\er2\}$. The
realized decompositions are listed in Table~\ref{hyperbolic}.
\endtheorem

As one can easily see, the Pontrjagin-Viro form of an $M$-surface of
hyperbolic type is uniquely recovered from~\eqref{E-RGM}. The
corresponding complex separations and values $\CP(w_1)$ are given in
Table~\ref{hyperbolic}.

\midinsert
\table\label{hyperbolic}
$M$-surfaces of hyperbolic type
\endtable
\vglue-\medskipamount
\let\style\textstyle
\def\gap{\quad}
\def\\{\noalign{\smallskip}}
\catcode`\+\active \let+\+
\let\0\varnothing
\def\c{0}
\def\n{2}
\def\cn{\hidewidth0,2\hidewidth}

\hbox to\hsize{\eightpoint\hss\vtop{%
\halign{\llap{$\style^#{}\,$}&$\style#$\hss&\gap$\style#$\hss\cr
\multispan3\hss\bf Case $E_{\R}=V_3+V_1+4S$\hss\cr
\noalign{\smallskip}
&(V_3\+V_1)\+(\0)    &(2S)\+(2S)    \cr
&(V_3\+S)\+(\0)      &(V_1\+S)\+(2S)\cr
&(V_3\+S)\+(V_1)     &(2S)\+(S)     \cr
&(V_3\+S)\+(S)       &(V_1\+S)\+(S) \cr
&(V_3\+V_1\+S)\+(S)  &(S)\+(S)      \cr
&(V_3\+2S)\+(S)      &(V_1)\+(S)    \cr
&(V_3\+2S)\+(V_1\+S) &(S)\+(\0)     \cr
&(V_3\+2S)\+(2S)     &(V_1)\+(\0)   \cr
&(V_3\+V_1\+2S)\+(2S)&\0            \cr}}%
\hss\vrule\hss\vtop{%
\halign{\llap{$\style^#{}\,$}&$\style#$\hss&
  \gap$\style#$\hss&\gap\hss$\style#$\hss\cr
\multispan4\hss\bf Other cases\hss\cr
\noalign{\smallskip}
&(V_4\+S)\+(\0)     &(2S)\+(2S) &\c\cr
\noalign{\smallskip}
&(V_{11}\+V_1)\+(\0)&\0            \cr
&(V_{11})\+(\0)     &(V_1)\+(\0)   \cr
\noalign{\smallskip}
&(V_{10})\+(\0)     &(V_2)\+(\0)&\c\cr
&(V_{9})\+(\0)      &(V_3)\+(\0)   \cr
&(V_{8})\+(\0)      &(V_4)\+(\0)&\n\cr
&(V_{7})\+(\0)      &(V_5)\+(\0)   \cr
&(V_{6})\+(\0)      &(V_6)\+(\0)&\c\cr
\noalign{\smallskip}
&(V_{10})\+(\0)     &(S_1)\+(\0)&\c\cr
}}\hss}
\endinsert

\subsection{Other surfaces with Pontrjagin-Viro form}\label{7.3}

\theorem\label{th-other}
Nonmaximal real Enriques surfaces admitting Pontrjagin-Viro form are
those and only those listed in Table~\ref{other}. The Pontrjagin-Viro
form of such a surface is determined, via~\eqref{E-RGM}, by the
decomposition $\er=\{\er1\}\+\{\er2\}$. Table~\ref{other} lists
the complex separations and values $\CP(w_1)$ on the characteristic
class of a nonorientable component of even genus, if such a component
is present. A~${}^*$ marks the decompositions which are also realized
by a real Enriques surface of type~$\II$\rom; a ${}^{**}$ marks the
decompositions which are also realized by a real Enriques surface of
type~$\I{}$ not admitting Pontrjagin-Viro form.
\endtheorem

\proof\nofrills{\/}
of Theorem~\ref{th-other} is based upon the classification of real
Enriques surfaces, which will appear in full in~\cite{DIK} (see also
\cite{DK1}--\cite{DK3}). The necessary partial results are cited
below.

\midinsert
\table\label{other}
Other surfaces with Pontrjagin-Viro form
\endtable
\vglue-\medskipamount
\let\style\textstyle
\def\gap{\quad}
\def\\{\noalign{\smallskip}}
\catcode`\+\active \let+\+
\let\0\varnothing
\def\c{0}
\def\n{2}
\def\cn{\hidewidth0,2\hidewidth}

\hbox to\hsize{\eightpoint\hss\vtop{%
\halign{\llap{$\style^#{}\,$}&$\style#$\hss&
  \gap$\style#$\hss&\gap\hss$\style#$\hss\cr
\multispan4\hss\bf$(M-2)$-surfaces\hss\cr
\noalign{\smallskip}
&(V_4)\+(\0)	&(V_1)\+(V_1)	&\c\cr
&(V_4)\+(V_1)	&(V_1)\+(\0)	&\c\cr
&(V_4)\+(2V_1)	&\0		&\c\cr
\noalign{\smallskip}
&(V_3)\+(V_1)	&(V_2)\+(\0)	&\n\cr
&(V_3)\+(\0)	&(V_2)\+(V_1)	&\n\cr
\noalign{\smallskip}
&(V_6)\+(\0)	&(S)\+(S)	&\c\cr
\noalign{\smallskip}
&(V_5)\+(\0)	&(V_1)\+(S)	   \cr
&(V_5)\+(V_1)	&(S)\+(\0)	   \cr
&(V_5)\+(S)	&(V_1)\+(\0)	   \cr
&(V_5\+V_1)\+(S)&\0		   \cr
\noalign{\smallskip}
&(V_4)\+(\0)	&(V_2)\+(S)	&\n\cr
&(V_4)\+(S)	&(V_2)\+(\0)	&\c\cr
\noalign{\smallskip}
&(V_3)\+(S)	&(V_3)\+(\0)	   \cr
\noalign{\smallskip}
&(V_4)\+(S)	&(S_1)\+(\0)	&\c\cr
\noalign{\smallskip}
  &(V_2)\+(V_2)		&(S_1)\+(\0)	&\n\cr
\crcr}}%
\hss\vrule\hss\vtop{%
\halign{\llap{$\style^{#}{}\,$}&$\style#$\hss&
  \gap$\style#$\hss&\gap\hss$\style#$\hss\cr
\multispan4\hss\bf$(M-2)$-surfaces (continued)\hss\cr
\noalign{\smallskip}
 *&(V_1)\+(V_1)		&(2S)\+(S)	\cr
 *&(2V_1)\+(S)		&(S)\+(S)	\cr
 *&(V_1\+S)\+(V_1\+S)	&(S)\+(\0)	\cr
 *&(2V_1\+S)\+(2S)	&\0		\cr
 *&(V_1\+S)\+(S)	&(V_1)\+(S)	\cr
 *&(V_1\+S)\+(2S)	&(V_1)\+(\0)	\cr
\noalign{\smallskip}
**&(V_2\+2S)\+(2S)	&\0		&\c\cr
  &(V_2\+S)\+(S)	&(S)\+(S)	&\c\cr
  &(V_2)\+(\0)		&(2S)\+(2S)	&\c\cr
\noalign{\smallskip}
**&(S_1)\+(\0)		&(2S)\+(2S)	\cr
\noalign{\smallskip}
**&(V_{10})\+(\0)	&\0		&\c\cr
\noalign{\medskip}
\multispan4\hss\bf$(M-4)$-surfaces\hss\cr
\noalign{\smallskip}
**&(V_4)\+(S)		&\0		&\c\cr
\noalign{\smallskip}
**&(2S)\+(2S)		&\0		\cr
  &(S)\+(S)		&(S)\+(S)	\cr
\crcr}}\hss}
\endinsert

The fact that in all cases listed in Table~\ref{other} the
Pontrjagin-Viro form is determined by~\eqref{E-RGM} is
straightforward. Thus, it remains to enumerate the surfaces for which
the Pontrjagin-Viro form is well defined. In view of~\eqref{A-congr}
for such a surface~$E$ one has $\chi(\er)=8$, $0$, or~$-8$. If
$b_0(\er)=1$, \ref{E-RGM} applied to an empty quoter gives
$\chi(\er)=-8$. Thus, it suffices to consider $(M-d)$-surfaces with
either $\chi(\er)=8$ and $d=2,4$, or $\chi(\er)=0$ and $d=2,4$, or
$\chi(\er)=-8$ and $d=2$.

\case1{$\chi(\er)=-8$, $d=2$.}
The only topological type $\er=V_{10}$. {\proclaimfont There are two
deformation families of real Enriques surfaces~$E$ with
$\er=V_{10}$\rom; they are both of type~$\I{}$ and differ by whether
$w_2(E/\conj)$ is or is not~$0$ \rm(see~\cite{DK3} and~\cite{DIK}).}
Since $\er$ has a single component, such a surface admits the
Pontrjagin-Viro form if and only if $w_2(E/\conj)=0$.

\case2{$\chi(\er)=0$, $d=2$.}
All such surfaces are of type~$\I{}$ (see~\cite{DK2}); hence, they
satisfy the hypotheses of \iref{sufficient}3 and the Pontrjagin-Viro
form is well defined.

\case3{$\chi(\er)=0$, $d=4$.} Among the five topological types, with
$\er=2S_1$, $S_1\+V_2$, $2V_2$, $V_3\+V_1$, and $V_4\+S$
(see~\cite{DK1}), only the last one can satisfy~\ref{E-RGM}. (Recall
that an $S_1$ component must form a separate half.) Furthermore, the
complex separation must be $\{(V_4)\+(S)\}\+\{\varnothing\}$; in
particular, both the components of~$\er$ are in one half.
{\proclaimfont There are two deformation families of real Enriques
surfaces~$E$ with $\er=\{V_4\+S\}\+\{\varnothing\}$
{\rm(see~\cite{DIK})}.\footnote{In~\cite{DK3} it is erroneously
stated that there is one family.} They can be obtained by the
construction of~\ref{5.3} from a triple $(\tZ;P,Q)$, where
$\tZ=\Sigma_4$ is a rational ruled surface with a $(-4)$-section,
$P\in|2e_\infty|$, and $Q\in|e_0+e_\infty|$.  \rom(Here $e_0$ is the
exceptional $(-4)$-section and $e_\infty$ is the class of a generic
section.\rom) One type is obtained when $\tZ_\R=S_1$ and
$P_\R=\varnothing$\rom; the other one, when $\tZ_\R=\varnothing$.}
In the former case, $\tZ_\R=S_1$, one can apply
Theorem~\ref{main.th}: since $P$ is of type~$\II$, the
Pontrjagin-Viro form is not defined. Thus, it suffices to construct a
surface with $\er=V_4\+S$ and well defined Pontrjagin-Viro form. This
can be done as in~\ref{5.3} and~\ref{5.7}, where
$\tZ'=\Rp1\times\Rp1$, $P'$ is a pair of conjugate generatrices of
bi-degree~$(0,1)$, and $Q'$ is the union of a pair of generatrices
of bi-degree~$(1,0)$ and a nonsingular curve of bi-degree~$(2,2)$.

\case4{$\chi(\er)=8$, $d=2$.} There are three topological types, with
$\er=2V_1\+3S$, $V_2\+4S$, and $S_1\+4S$ (see~\cite{DK1}), which we
consider separately.

{\proclaimfont Each decomposition $\er=\{\er1\}\+\{\er2\}$ of
$\er=2V_1\+3S$ is realized by two deformation families of real
Enriques surfaces, one of type~$\I{}$, and one of type~$\II$
\rm(see~\cite{DIK}).} The surfaces of type~$\I{}$ satisfy the
hypotheses of~\iref{sufficient}3 and, hence, have well defined
Pontrjagin-Viro forms.

{\proclaimfont There is one deformation family of real Enriques
surfaces~$E$ with $\er=\{V_2\}\+\{4S\}$, one family of surfaces with
$\er=\{V_2\+2S\}\+\{2S\}$, and two families of surfaces with
$\er=\{V_2\+4S\}\+\{\varnothing\}$ {\rm(see~\cite{DK3}
and~\cite{DIK}).\footnote{In~\cite{DK3} it is erroneously stated that
there is one family with $\er=\{V_2\+4S\}\+\{\varnothing\}$.}} All
these surfaces are of type~$\I{}$.} In the first two cases the
surfaces satisfy the hypotheses of~\iref{sufficient}3 and, hence,
have well defined Pontrjagin-Viro forms. In the last case the
surfaces can be obtained by the construction of~\ref{5.3}
and~\ref{5.7}:  one takes for~$\tZ'$ the projective plane~$\Rp2$,
for~$P'$, the union of two conics with two conjugate tangency points
(so that $P_\R=\varnothing$), and for~$Q'$, the union of a generic
real line and the line through the singular points of~$P'$. The
conics of~$Q'$ may be either both real or complex conjugate; in the
former case $Q'$ is of of type~$\II$ and the Pontrjagin-Viro form is
not defined, in the latter case $Q'$ is of type~$\I{}$ and the
Pontrjagin-Viro form is defined due to Theorem~\ref{main.th}.

{\proclaimfont There are two deformation families of real Enriques
surfaces~$E$ with $\er=\{S_1\}\+\{4S\}$ {\rm(see~\cite{DK3}
and~\cite{DIK}).} They are obtained by the construction of~\ref{5.3}
and~\ref{5.7} from a triple $(\tZ';P',Q')$, where $\tZ'$ is the plane
$\Rp2$ \rom(or hyperboloid $\Rp1\times\Rp1$\rom), $P'$ is a
nonsingular $M$-curve of degree~$4$ \rom(respectively,
bi-degree~$(4,2)$\,\rom), and $Q'$ is a pair of conjugate lines
\rom(respectively, generatrices of bi-degree~$(0,1)$\,\rom).} From
Theorem~\ref{main.th} it follows that the Pontrjagin-Viro form is
well defined in the latter case and is not defined in the former
case (as $\tZ^-\sminus Q_\R$ is nonorientable).

\case5{$\chi(\er)=8$, $d=4$, i.e., $\er=4S$.}
Only the decompositions $\{4S\}\+\{\varnothing\}$ and
$\{2S\}\+\{2S\}$ can satisfy~\ref{E-RGM} (and, in fact, only
these surfaces are of type~$\I{}$). Consider the two cases
separately.

{\proclaimfont There are four deformation families of real Enriques
surfaces~$E$ with $\er=\{4S\}\+\{\varnothing\}$
{\rm(see~\cite{DIK})}. They differ by the classes realized by the
image of~$\xr1$ in $X/\tau=E$ and $X/\c2$
\rom(see~\ref{5.1}\rom)\rom; the four possibilities are $(w_2,w_2)$,
$(w_2,0)$, $(0,w_2)$, and $(0,0)$.} (Note that $X/\c2$ is
diffeomorphic to an Enriques surface and $w_2(X/\c2)\ne0$.)
Since $E/\conj$ can as well be represented as the quotient space of
$X/\c2$ by an involution whose fixed point set is $\xr1/\c2$, only
the first of the four families may possess Pontrjagin-Viro form. Such
a surface can be obtained by the construction of~\ref{5.3}
and~\ref{5.7}. Take for~$\tZ'$ the hyperboloid $\Rp1\times\Rp1$. Let
$(L_1,L_2)$ and $(M_1,M_2)$ be two pairs of conjugate generatrices of
bi-degree~$(1,0)$ and $(N_1,N_2)$ a pair of conjugate generatrices of
bi-degree~$(0,1)$. Pick a generic pair~$(C_1,C_2)$ of conjugate
members of the pencil generated by $L_1+M_1+N_1$ and $L_2+M_2+N_2$
and let $P'=C_1+C_2$ and $Q'=N_1+N_2$. ($P'$ is a real curve of
type~$\I{}$ with four nodes, which lie on~$Q'$.) Existence of the
Pontrjagin-Viro form follows from Theorem~\ref{main.th}.

{\proclaimfont There is one deformation family of real Enriques
surfaces~$E$ with $\er=\{2S\}\+\{2S\}$ \rm(see~\cite{DIK}).} A
surface with Pontrjagin-Viro form is constructed similar to the
previous case. One takes for~$Z'$ the hyperboloid $\Rp1\times\Rp1$,
for~$P'$, a real $M$-curve of bi-degree~$(4,2)$ with two conjugate
double points, and for~$Q'$, the union of two conjugate generatrices
through the singular points of~$P'$.
\endproof


\widestnumber\key{DK2}
\Refs

\ref{vdB}\label{vdBlij}
\by	van der Blij
\paper	An invariant of quadratic forms $\!{}\bmod8$
\jour	Indag. Math.
\vol	21
\yr	1959
\pages	291--293
\endref

\ref{Br}\label{Brown}
\by 	E.~H.~Brown
\paper	Generalization of the Kervaire invariant
\jour	Ann. Math.
\vol	95
\yr	1972
\pages	368--383
\endref

\ref{CM}\label{CM}
\by P.~E.~Conner, E.~Y.~Miller
\book Equivariant self-intersection
\bookinfo Preprint
\yr 	1979
\endref

\ref{GM}\label{GM}
\by L.~Guillou, A.~Marin
\paper Une extension d'un th\'eor\`eme de Rokhlin
sur la signature
\jour C.~R.~Acad. Sci. Paris S\'er.~A
\vol 	285
\yr 	1977
\pages	95--98
\endref

\ref{DIK}\label{DIK}
\by A. Degtyarev, I. Itenberg, V. Kharlamov
\book Real Enriques surfaces
\toappear
\endref

\ref{DK1}\label{DK1}
\by A.~Degtyarev, V.~Kharlamov
\paper Topological classification of real Enriques surfaces
\jour Topology
\vol   35
\yr    1996
\issue 3
\pages 711--729
\endref

\ref{DK2}\label{DK2}
%\by A.~Degtyarev, V.~Kharlamov
\bysame
\paper	Halves of a Real Enriques Surface
\jour	Comm. Math. Helv.
\vol	71
\yr	1996
\pages	628--663
\moreref\nofrills Extended version:
\book Distribution of the components of a real Enriques
surface
\bookinfo Preprint of the Max-Planck Institute, MPI/95-58
\yr  1995
\miscnote also available from {\smc AMS} server as
{\smc AMSPPS} \#~199507-14-005
\endref

\ref{DK3}\label{DK3}
%\by A.~Degtyarev, V.~Kharlamov
\bysame
\paper	On the moduli space of real Enriques surfaces
\jour	C.~R.~Acad. Sci. Paris S\'er~\rom{I}
\yr	1997
\endref

\ref{EM}\label{ME}\label{Pontrjagin}
\book	Encyclopedia of Mathematics
\publ	Reidel Kluver Academic Publishers
\publaddr P.O.Box 17 3300 1A Dordreecht Holland
\yr	1988
\endref

\ref{Fin}\label{Finashin}
\by S.~Finashin
\paper $\text{Pin}^-$-cobordism invariant and generalization of
Rokhlin signature congruence
\jour Algebra i Analiz
\vol 	2
\yr 	1990
\issue	2
\pages	242--250
\lang Russian
\transl\nofrills English transl. in
\jour Leningrad Math. J.
\vol 	2
\yr 	1991
\issue	2
\pages	917--924
\endref

\ref{Ka}\label{Kalinin}
\by I.~Kalinin
\paper Cohomological characteristics of real projective
hypersurfaces
\jour  Algebra i Analiz
\vol 	 3
\yr 	  1991
\issue 2
\pages 91--110
\lang Russian
\moreref\nofrills English transl. in
\jour  St.~Petersburg Math.~J.
\vol   3
\yr    1992
\issue 2
\pages 313--332
\endref

\ref{K\"u}\label{Kucuk}
\by \"O.~K\"u\cedilla c\"uk
\book
\bookinfo M.S.~thesis
\miscnote in preparation
\endref

\ref{KV}\label{KV}
\by	V.~Kharlamov, O.~Viro
\paper	Extensions of the Gudkov-Rokhlin congruence
\jour	Lect. Notes in Math.
\vol	1346
\yr	1988
\pages	357--406
\endref

\ref{Mik}\label{Mikhalkin}
\by G.~Mikhalkin
\paper The complex separation of real surfaces and extensions of
Rokhlin congruence
\jour	Invent. Math.
\vol	118
\yr	1994
\pages	197--222
\endref

\endRefs
\enddocument